\documentclass[3p]{elsarticle} 

\usepackage{graphicx}
\usepackage{subfig}
\usepackage{array} 
\usepackage{amssymb}
\usepackage{amsmath,mathtools}
\usepackage{bm}
\usepackage{float}
\usepackage{calc}
\usepackage{color}
\usepackage{url}
\usepackage{ulem}
\usepackage{soul}
\usepackage{multirow}
\usepackage{booktabs}
\usepackage{comment}

\newtheorem{theorem}{\textbf{Theorem}}[section]

\newtheorem{proposition}[theorem]{Proposition}

\makeatletter
\newcommand\avsuminner[2]{%
  {\sbox0{$\m@th#1\sum$}%
   \vphantom{\usebox0}%
   \ooalign{%
     \hidewidth
     \smash{\vrule height\dimexpr\ht0+1pt\relax depth\dimexpr\dp0+1pt\relax}%
     \hidewidth\cr
     $\m@th#1\sum$\cr
   }%
  }%
}
\makeatother

\begin{document}

\title{A Reactive Power Market for the Future Grid}
\tnotetext[t1]{This material is based upon work partially supported by the MIT Energy Initiative, MathWorks Mechanical Engineering Fellowship, Chevron-MIT Energy Fellowship, and the U.S. Department of Energy under Award Number DE-IA0000025. The views and opinions of authors expressed herein do not necessarily state or reflect those of the United States Government or any agency thereof.}

\author[1]{Adam Potter}
\ead{awpotter@mit.edu}

\author[1]{Rabab Haider\corref{cor1}}
\ead{rhaider@mit.edu}

\author[2]{Giulio Ferro}
\ead{giulio.ferro@unige.it}

\author[2]{Michela Robba}
\ead{michela.robba@unige.it}

\author[1]{Anuradha M. Annaswamy}
\ead{aanna@mit.edu}
	
\cortext[cor1]{Corresponding author}

\address[1]{Department of Mechanical Engineering, Massachusetts Institute of Technology, Cambridge, MA, USA}
\address[2]{Department of Informatics, Bioengineering, Robotics and Systems Engineering, University of Genoa, Genoa, Italy}

\begin{abstract}
As pressures to decarbonize the electricity grid increase, the grid edge is witnessing a rapid adoption of distributed and renewable generation. As a result, traditional methods for reactive power management and compensation may become ineffective. Current state-of-art for reactive power compensation, which rely primarily on capacity payments, exclude distributed generation (DG). We propose an alternative: a reactive power market at the distribution level designed to meet the needs of decentralized and decarbonized grids. The proposed market uses variable payments to compensate DGs equipped with smart inverters, at an increased spatial and temporal granularity, through a distribution-level Locational Marginal Price (d-LMP). We validate our proposed market with a case study of the US New England grid on a modified IEEE-123 bus, while varying DG penetration from 5\% to 160\%. Results show that our market can accommodate such a large penetration, with stable reactive power revenue streams. The market can leverage the considerable flexibility afforded by inverter-based resources to meet over 40\% of reactive power load when operating in a power factor range of 0.6 to 1.0. DGs participating in the market can earn up to 11\% of their total revenue from reactive power payments. Finally, the corresponding daily d-LMPs determined from the proposed market were observed to exhibit limited volatility.

\end{abstract}


\maketitle

\section{Introduction} \label{sec:intro}
 
With a growing demand for zero-carbon generation in the electricity sector comes an increase in grid-edge distributed energy resources (DERs)\footnote{\textit{Abbreviations:}
CI - Current Injection, CI-OPF - Current Injection optimal power flow, d-LMP -, Distribution Locational Marginal Price, DER - Distributed Energy Resource, DG - Distributed Generator, DSO - Distribution System Operator , ISO Independent System Operator, LMP - Locational Marginal Price, LOC - Lost Opportunity Cost, MCE - McCormick Envelope, OPF - Optimal Power Flow, P - Real Power, PAC - Proximal Atomic Coordination , PF - Power Factor, PMO - Primary Market Operator, PV - Photovoltaic, Q - Reactive Power, RTO - Regional Transmission Organization}, and in particular, a subset of DERs referred to as distributed generation (DGs). The increasing penetration of these small-scale resources disrupts the existing industry practice for power dispatch, device control, and market compensation mechanisms. Traditionally, large controllable generators at the transmission level have been responsible for maintaining grid stability and power quality; however new strategies are required for the emerging grid. While recent efforts such as FERC Order 2222 \cite{ferc_2222} in the US indicate a move towards accommodation of DERs, new electricity market structures are needed to fully utilize the flexibility they offer at high penetration and lead to a feasible framework for their integration while preserving grid stability and power quality. This paper identifies the requirements for such a market structure, and proposes one with a focus on a reactive power market at the distribution level.
 
 \begin{figure}[h!]
    \centering
    \includegraphics[scale=0.5]{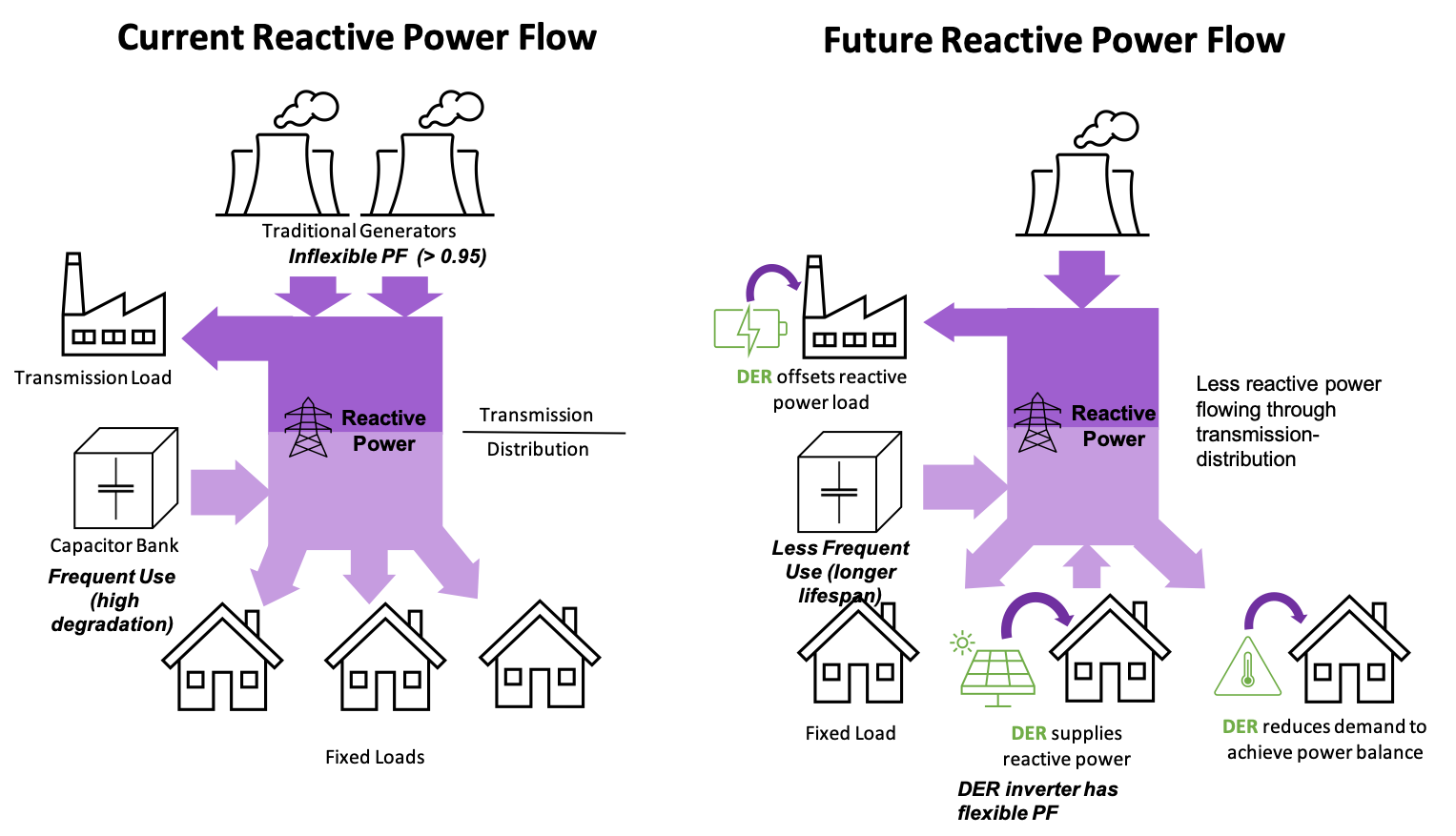}
    \caption{Balancing reactive power in an electric grid is necessary to stabilize voltage and ensure power transfer. As grids decentralize, DERs can be utilized to provide reactive power support throughout the distribution system, with more flexible power factors and distributed control units. They can replace traditional reactive power suppliers, which are limited in PF and expensive to operate. However, to provide reactive power DERs must be incentivized to do so.}
    \label{fig:q_flow}
\end{figure}
 
Reactive power is intimately connected with power quality, and its regulation is achieved in the transmission grid by altering the power factors (PFs) of large synchronous generators. The outputs of these generators are easily controllable thus permitting an efficient regulation of PF, and have therefore been sufficient for ensuring power quality and voltage stability at the transmission level. PF can be regulated at the requisite level even with a few large synchronous generators  as the target PF value rarely drops below 0.95 \cite{zhong_reactive_2002}. In contrast to this practice at the transmission level, power quality in the distribution grid is typically achieved through devices such as on-load tap changers, capacitor banks, and voltage regulators, which locally inject reactive power. Such a solution however becomes quickly inadequate in a distribution grid with an increased penetration of DERs, which induce new consumption patterns and reverse power flows at faster timescales than previously experienced. Emerging solutions that allow regulation in a wide range of PF correspond to inverter-based resources, such as DGs which are often equipped with smart inverters, and can adjust their PF using suitable power electronics control, quickly and at low costs. Further, since reactive power travels poorly, the spatially distributed support offered by these DERs offers a unique grid flexibility option if their reactive power capabilities can be coordinated. This paradigm shift in reactive power support from centralized large-scale traditional generators to distributed small-scale consumer owned DERs is illustrated in Fig.~\ref{fig:q_flow}.

This operational change implies that DGs that enable such an additional flexibility must be compensated for through an appropriate market mechanism. A market structure that compensates DGs for their services, particularly for reactive power, is structurally different from the current practice \cite{federal_energy_regulatory_commission_payment_2014}. The development of such a market is the topic of this paper, wherein we propose an energy market for the distribution grid which prices reactive power. We build on the authors' earlier work in \cite{HaiderRM, Haider_ADAPEN} which proposed an energy market that enables fine granularity, in both location and in time, of power injections and corresponding energy prices. Unlike those works, here we focus primarily on reactive power, including the deployment of DER resources to provide reactive power, the financial remuneration for the corresponding grid services, and the role of the market in an unbalanced distribution grid with increasing DER penetration. We emphasize that our proposal is for an \textit{energy market} which \textit{simultaneously} prices both real and reactive power. The design of additional market derivatives, such as pricing for reactive power support in voltage control applications, remain outside the scope of this paper as ancillary services. A review of voltage control strategies can be found in \cite{mahmud_review_2016}.

Currently, reactive power markets do not exist in the US -- neither energy or ancillary markets. In Europe, approaches to ancillary services vary but none are inclusive to DGs despite recent efforts to reduce minimum power bids from 10 MW to 1 MW in several countries \cite{rancilio_ancillary_2022}.  For both regions, and many others, transitioning ancillary service structure to be inclusive of DGs is a proven priority for regulatory agents. In a common approach to reactive power compensation, and the approach used in the US, any payments to generators that provide reactive power are static rates, designed to recover the capital costs of installing the equipment or based on the LOC (lost opportunity cost) of not selling real power. These payments can consist of capacity, or fixed payments, based on generator's reactive power capabilities; and utilization, or variable, payments based on the delivered kVARs of reactive power. These static rates do not reflect the temporal shifts in supply and demand, or spatial variation across the grid, both of which are pertinent characteristics of DGs. Further, these payments are typically cleared on a case-by-case basis for each generator, a process which is generally exclusive to larger generators. This approach of static case-by-case rates will be insufficient in accommodating a deep penetration of DGs. These DGs offer significant operational flexibility at little to no cost, and can be leveraged to provide much needed grid services through the clean energy transition. To do so, a market mechanism which prices their reactive power capabilities and their spatial and temporal variations is vital. This is possible with approaches based on optimal power flow (OPF), which allow the simultaneous determination of prices corresponding to both real and reactive power injections, enables an increasing penetration of DERs, and accommodates grid physics. In this paper, we propose a reactive power market that allows variable payments based on an OPF approach.


The reactive power market proposed in this paper is a part of a distribution-level market with oversight assumed to be provided by a distribution system operator (DSO). A DSO-centric market has several advantages, the first of which is that it enables DERs to directly participate in a local setting, and have the market location be collocated electrically, such as at a primary feeder node or a secondary feeder node. It also avoids tier bypassing, which can arise if DERs are directly dispatched by wholesale markets in response to the bulk grid objectives without considering local constraints \cite{kristov_tale_2016}. Our proposed reactive power market sits at the primary feeder level (4 to 35 kV level), with any resources and uncontrollable loads residing in the secondary feeders and below aggregated upward. Each Primary Market Operator (PMO) is assumed to have oversight over one primary feeder, and determines the schedule for DERs in terms of both real and reactive power injections as well as the retail electricity price of both real and reactive power, which we denote as a distribution-level Locational Marginal Price (d-LMP). Each PMO pays the participating DGs in that feeder at the d-LMP. The d-LMP also corresponds to the payment made by the loads to the PMO. As we effectively split the grid into smaller service regions, the proposed retail market has the potential to lead to less expensive and more extensive procurement and payment for reactive power resources. The study in \cite{zhong_localized_2004} corroborates this statement as well. A DSO can be tasked with the oversight of these PMOs and perhaps be commissioned as a non-profit entity, similar to ISO/RTOs, by the corresponding Public Utility Commission.

The underlying framework utilized to determine the market schedules and prices is based on constrained optimization, an initial version of which was proposed in \cite{HaiderRM,Haider_ADAPEN}. In the present work, we extend \cite{HaiderRM,Haider_ADAPEN} in two ways: (1) we consider a multi-phase unbalanced distribution grid; and (2) we carry out a detailed treatment of reactive power deployment and pricing. The case study in \cite{HaiderRM,Haider_ADAPEN} considered a specific level of penetration of DG of approximately 70\%, consisting entirely of real power generation. In this paper, we significantly extend the scope, by modeling each DG as an inverter-based resource with reactive power capabilities and variable PF, over a range of DG penetration from 5\% to 160\%. We do not consider reactive power support from flexible consumption or storage devices, though our market structure can be extended to accommodate them as well. We also consider an unbalanced distribution grid, not considered in \cite{HaiderRM,Haider_ADAPEN}, using a current injection (CI) model, and include analytical support for the reactive power price resulting from this model. We place our reactive power market in context of industry needs. Through numerical case study, we demonstrate how our market effectively utilizes DERs to provide appropriate local reactive power support while optimizing network-level objectives subject to grid constraints. We validate the market's flexibility and performance over a range of conditions, including DG penetration and inverter PF and discuss the implications of our market in the key areas of price volatility and supporting investment decisions.

The following are the specific contributions of the paper.

\begin{itemize}
\item We present an analysis of the state-of-art of reactive power compensation in theory and in practice. We outline how DERs disrupt traditional reactive power management, warranting new approaches to reactive power dispatch and pricing. We identify three priorities for a DER-inclusive reactive power market:
\begin{itemize}
    \item Utilize the significant reactive power capabilities of a large penetration of DERs and their wider PF operating range 
    \item Support efficient grid operation by minimizing grid losses using local DERs
    \item Provide fair remuneration to DERs for their grid services
    \item Incentivize long-term growth of DER adoption through market inclusion and corresponding reliable revenue streams
\end{itemize}
These constitute the motivation for our reactive power market.

\item We propose a reactive power energy market for DGs in the distribution grid which meet the priorities outlined above. The market dispatch and price clearing are solutions of a linear CI model of the unbalanced distribution grid. The market decisions are arrived at using a Proximal Atomic Coordination (PAC) based distributed optimization algorithm, demonstrated to converge to the optimal point under suitable conditions in \cite{RomvaryTAC,HaiderRM}. We further analyze the reactive power price resulting from the CI model, building the case for minimizing line losses as a network-wide view of reactive power.
\item Through a numerical case study of the
US New England grid modeled as a modified IEEE-123 bus, we demonstrate the following features:

\begin{itemize}
    \item \textit{Our market can accommodate a large amount of DGs:} We illustrate how reactive power support can be realized even with a high penetration of DGs. We show in particular that the revenue stream for reactive support to these DGs remain steady even as DG penetration increases from 5\% to 160\%, which indicates the feasibility of the proposed market at high penetrations, and even at low penetrations.
    \item \textit{Our market uses all available DGs efficiently:} We utilize the full flexibility of DGs to operate over a range of PF, from 0.6 to 0.95, to meet over 45\% of reactive power load, even with the reactive power price per unit kVAR remaining the same. The corresponding increase in the requisite reactive power injection is appropriately compensated by an increasing percentage of the revenue from reactive power support compared to total revenue from DGs.
    \item \textit{Our market supports grid voltages:} The reactive power support from DGs spatially distributed throughout the grid enables grid voltages to remain within allowable operating range even at high DG penetration. Further, the reactive power support enhances grid operation by increasing the mean voltage across the network and providing necessary support during periods of high demand.
    \item \textit{Our market is stable and provides reliable revenue:} A detailed assessment of the price variations across all DGs in the IEEE 123-bus over a week shows that the price variations primarily follow load fluctuations and are otherwise minimal.
\end{itemize}

\end{itemize}

The remaining paper is organized as follows. Section \ref{sec:review} presents an analysis of the state-of-art in both theory and practice, providing a historical perspective on our identified priorities for reactive power compensation in the US. Section \ref{sec:market-design} introduces our OPF-based reactive power pricing. Section \ref{sec:method} outlines the simulation setup and methods. Section \ref{sec:results} presents results of numerical evaluation of the proposed reactive power market on an IEEE-123 bus network in the context of our market design priorities. Section \ref{sec:conclusion} provides concluding remarks.



\section{State of Art of Reactive Power Compensation} \label{sec:review}
In the following sections we provide an analysis of the state-of-art for reactive power management and pricing in the US\footnote{A review of reactive power markets in EU markets is presented in \cite{Wolgast_reactivereview}}. We present an analysis of both industry practices and technical solutions. We identify the design objectives for reactive power pricing for traditional power systems, and elucidate how they are inadequate for the emerging grid with high DER penetration. Finally, we identify the priorities for the emerging grid, which our proposed market in Section~\ref{sec:market-design} is shown to satisfy. A summary of the analysis is presented in Table~\ref{tab:state_of_art_summary}.

\begin{table*}[h]
\begin{tabular}{>{\raggedright\arraybackslash}p{0.18\linewidth}>{\raggedright\arraybackslash}p{0.25\linewidth}|>{\raggedright\arraybackslash}p{0.18\linewidth}>{\raggedright\arraybackslash}p{0.25\linewidth}}
\multicolumn{2}{>{}c|}{Current Practice (State-of-Art)} & \multicolumn{2}{>{}c}{Proposed Practice for Emerging Grid (Proposed)} \\ \toprule
Methodology for Q pricing & Limitations of state-of-art & Methodology for Q pricing & Features of proposed method \\ \midrule
 Fixed pricing & Emerging grid is more dynamic: fixed prices cannot provide adequate signals for devices to respond to & Dynamic pricing &  Time-varying and spatially varying Q prices provide signal for DERs to respond to system needs \\
 Cost-based pricing & Not inclusive to DERs which don't incur high costs for producing Q  & Service-based pricing & Q price is determined by an OPF, thereby modeling the real-time value of Q injection to the grid \\
 Case-by-case clearing & Intractable for a large number of devices (i.e. DER-rich grid) & Coordinated and distributed clearing & Maintain tractability of large-scale OPF by using distributed optimization \\
 Price deadbands where no payment is made (0.95 to 1 PF) & Deadbands are anti-competitive and not recommended for inclusive markets  & No price deadbands & Provide compensation through the full range of PFs enabled by power electronic control \\ 
\end{tabular}
\caption{Summary of the state-of-art industry practice on reactive power management and pricing. We propose a reactive power market for the emerging grid, with new methodologies and features.}
\label{tab:state_of_art_summary}
\end{table*}

\subsection{History of reactive power pricing in the US}
Since the US energy sector was deregulated in the 1990s, there has been much debate over the treatment of reactive power in a market setting. Early proposals for compensating reactive power (Q) aimed to address two problems: 
\begin{itemize}
    \item [1.] Ensuring sufficient Q resources were available in the market: by focusing on the recovery of Q-enabling equipment costs;
    \item [2.] Maintaining grid stability: by efficiently procuring the requisite reactive power capacity in the long-term
\end{itemize}

These problems were inherent to the then centralized grid predominantly consisting of synchronous generators. Early reactive power compensation methods proposed in \cite{zhong_reactive_2002} and \cite{bhattacharya_reactive_2001} supplied no payment for generators operating at PFs in the range of 0.95 leading to 0.95 lagging creating a price ``deadband'', citing that synchronous generators do not incur losses within this operating range. Reactive power prices were set outside this PF range as it necessitated additional hardware and therefore recovery of equipment cost. Any further deviations imply a loss of opportunity of selling real power, and therefore the price is set based on the lost opportunity cost (LOC) \cite{rueda-medina_distributed_2013}. Much of current practice has followed this procedure. However, this approach is highly sensitive to the PF range and becomes inadequate as PF varies. Thus, this approach cannot be used for the emerging grid, where the penetration of DERs with variable PFs is rapidly increasing. Further, as clearly mentioned in \cite{federal_energy_regulatory_commission_payment_2014}, pricing deadbands for reactive power are anti-competitive and not recommended for inclusive markets, such as those encouraged in FERC Order 2222 \cite{ferc_2222}. On the demand side, some utilities enforce a similarly-structured PF penalty to customers that operate outside of a permitted PF range to recover costs from the additional electrical line losses incurred. However, this approach also uses a pricing deadband, is ineffective at signaling demand response actions, and can result in inequities between customers \cite{baughman_real-time_1991}.

In addition to efforts for recovering equipment costs and lost revenue, much of industry's attention was focused on the second problem of reactive power capacity procurement. As utilities rushed to ensure sufficient capacity to maintain voltage and power quality under high demand, capacity markets rose as a reliable way to incentivize reactive power capability investments. These fixed payments were preferred over variable payments for the VAR-hrs produced. References \cite{weber_simulation_1998} and \cite{barquin_gil_reactive_2000} use OPF-based methods with the underlying dual variables being used to generate reactive power prices, but both noted that the resulting prices were highly volatile. Study \cite{barquin_gil_reactive_2000} ultimately recommended capacity payments to drive investment. Reference \cite{ahmed_method_2000} assessed that before 1998 in the US around 80\% of reactive power payout was through fixed, or capacity, payments while only 20\% was through variable payments. It should be noted that even in this paper, which is more than 20 years old, fixed payments are recommended to be a transitional process phasing into payment purely on the utilization of reactive energy. For more than 20 years, researchers have argued fixed payments are outdated and should not remain the dominant compensation method for Q-based services. These fixed payments do not reflect the true value of services provided. Further, while synchronous generators are limited in their ability to generate Q, which can result in a capacity shortage, modern DERs such as PV solar equipped with smart inverters are not. As the penetration of inverter-based resources capable of generating a surplus of Q increases, the recruitment of Q-generating devices will no longer be a primary concern. Rather, the challenge of maintaining grid stability in an inverter-rich grid comes down to effectively utilizing the large number of available resources, and designing the corresponding payment structures required to incentivize Q-generation. 
Finally, and of critical importance, capacity markets are not inclusive to DERs. Therefore, as we move towards a grid with a strong influx of DG resources with much wider PF operating range, such procurement focused capacity markets will become less important. Like real power, it is crucial that reactive power compensation transitions to market-based solutions, and accounts for value added to the grid, the spatial differences in resources, and are inclusive to DERs. In this paper, we propose a specific distributed optimization based solution that allows the determination of a reactive power pricing structure that captures a fine-grain variation of location and time.

\subsection{Implementation practices for reactive power pricing}
In this section we focus on the Q-pricing strategies currently implemented in practice. Reactive power compensation at the transmission level is implemented using various methods. Since reactive power capabilities are deemed critical by FERC, monetary incentives are required to ensure power is delivered and the grid has sufficient capacity. Tariffs are designed as either fixed or variable payments to asset owners. However, poorly designed incentives could result in either a reactive power shortage resulting in voltage control failures, or surpluses that create monetary losses for utilities and consumers. A snapshot of these methods at various ISO/RTOs is delineated in this section.

A US FERC report released in 2014 outlines two common payment methods \cite{federal_energy_regulatory_commission_payment_2014}. The first is a regulatory mandate on generator reactive power support. Under the assumption that reactive power is inexpensive, all market participants are expected to help balance the grid without additional payment. This method is clearly outdated and unfair, and further, is incapable of incorporating DERs. The second method, developed by the American Electric Power Service (AEP), permits generators to report expenses explicitly tied to reactive power capabilities and settle on an appropriate fixed payment. The settlement process also gives grid operators control over capacity by recruiting and settling as necessary. This cost-based compensation is commonly used today, and is appealing because it avoids a market and the associated risks of market failures while ensuring the grid has appropriate capacity of reactive power. However, this method is generally limited to large traditional generators who can afford the cost analyses and settlement process. Additionally, unlike traditional generators, most DGs come equipped with extensive reactive power capabilities so expense-based pricing is less applicable. As mentioned earlier, the report \cite{federal_energy_regulatory_commission_payment_2014} alludes to the need for retail level accommodation of DERs as their penetration increases. Without a market, recruiting many small DERs for reactive power support with the industry standard would be prohibitive and the resulting prices unsubstantiated. A well designed reactive power market could facilitate the rapid adoption of small-scale distributed grid assets. 



\begin{figure}[h!]
    \centering
    \includegraphics[scale=0.45]{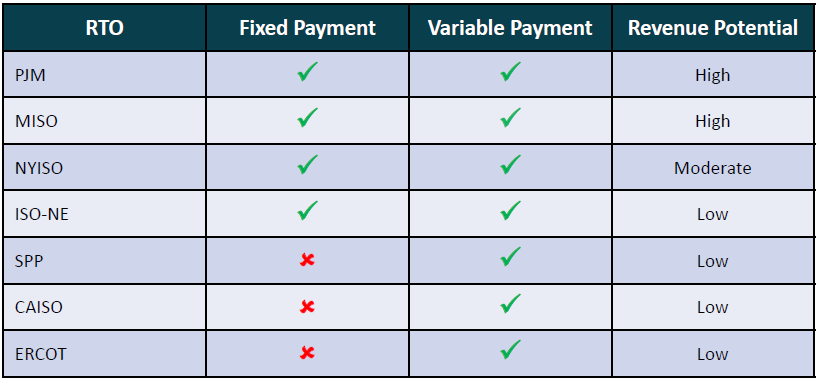}
    \caption{Reactive power compensation for ISO/RTO regions in the US with corresponding estimates on the reactive power revenue potential for solar and storage in those regions \cite{michael_borgatti_adrian_kimbrough_steven_shparber_reactive_nodate}. }
    \label{fig:revenue_potential}
\end{figure}

Another point of interest is the 2020 primer from the Solar Energy Industries Association (SEIA) \cite{michael_borgatti_adrian_kimbrough_steven_shparber_reactive_nodate}. The primer suggests variable payments are currently inadequate for incentivizing reactive power generation; in practice the ISO/RTOs rely heavily on fixed payments. Figure~\ref{fig:revenue_potential} outlines reactive power revenue potential for large solar projects. All US ISO/RTOs include variable, per kVAR-yr, payment based on loss of potential revenue in the real power market, but only regions with fixed payments provide high or moderate potential for a reactive power revenue stream. Currently, if a market participant, for example DERs, did not have access to fixed payments their revenue potential would be compromised.

The issue of loss of potential revenue in the real power market is further substantiated by data from the New England Independent System Operator (ISO-NE) in the US settled cost for yearly reactive power capacity. In 2020, the costs varied from 1.095 \$/kVAR-yr to 1.13 \$/kVAR-yr; in comparison, the normal utilization price is 0.05 \$/MWh (0.0057 \$/kVAR-yr). As the SEIA primer suggested, current tariff structures heavily favor standby capacity over power delivered. The report further notes that while the AEP method for reactive power capacity pricing has been used for several large solar projects, it is not inclusive to residential solar. Looking at ISO-NE again, the smallest asset reported had 587 kVARs of capacity, confirming the standard is non-inclusive to smaller DGs \cite{iso_express}. This implies that despite fixed payments being the go-to procedure in industry practice, there is still a significant gap in the efficient integration of small-scale resources. 
One pilot program running since 2017 in Italy aimed to test the feasibility of DER provision of ancillary services by opening the Italian ancillary service market to smaller DGs and DG aggregators \cite{gulotta_opening_2020}. The pricing approach had agents bid for yearly capacity payments that require DGs to bid for variable payments daily. While the market was overall more competitive, many DGs help capacity contracts but did not contribute much to the ancillary services as their variable price bids were too high. Here, capacity payments recruited DGs to participate, but the variable pricing scheme did not encourage their contribution to the grid.
A pricing scheme centered around fixed payments would be difficult to scale, improperly value resources, and result in tier bypassing by separating real and reactive power. 

The above analysis suggests a need for market-driven variable payments rather than fixed payments, to accommodate the services provided by DGs. In this paper, we propose a reactive power market that allows variable payments for resources in the distribution grid for energy services provided at the distribution level, and argue it is effective and reliable enough to facilitate the transition to high DER penetration.

\subsection{Technical Solutions for Reactive-power Pricing}\label{sec:review:lit}
A wide variety of pricing and settlement methods for reactive power have been proposed in the literature. The survey in \cite{Wolgast_reactivereview} reviews various approaches such as optimization, machine learning, and game theory, and presents a discussion of different market structures including ancillary provisions. Another review paper compares approaches to ancillary services emphasizing distinctions in the relationship between real and reactive power, the use of LOC and capacity pricing, the grid model used, and the spatial and temporal pricing granularity \cite{jay_comprehensive_2021}. In this paper, we focus on energy markets, and leave the development of ancillary markets at the distribution level - such as those for voltage control applications - to future work. The survey in \cite{safari_reactive_2020} presents a review of both LOC-based and OPF-based methods for reactive power pricing schemes. These can loosely be distinguished into three categories as described in \cite{homaee_retail_2020} based on the relationship between real and reactive power pricing. The first category includes the cost-based methods mentioned above that keep reactive power as an independent ancillary service. The dynamics of real and reactive power are assumed to be decoupled allowing voltage control to be a secondary concern to that of balancing power. Control of voltage is assumed to be possible through regulation of PFs using large generators. This solution is predicated on transmission-level resources, and the assumptions of decoupled real and reactive power are violated when PFs dip below 0.95, a plausible occurrence for DERs. This is further exacerbated by the unbalanced physics of the distribution grid.

The second category proposed in \cite{homaee_retail_2020} includes a hierarchical approach where the first step is the clearing of a real power market followed by a reactive power market. These solutions leverage the advantages of decoupled dynamics leading to simple market-based approaches. Reactive power is priced competitively and supports existing real power markets. One such example can be found in \cite{hasanpour_new_2009} which is based on LOC of selling real power, while others such as \cite{sarmila_har_beagam_market_2018} formulate a separate optimal reactive power flow problem to quantify its marginal value. In both cases, the use of a hierarchical structure may prevent resources that participated in the real power market from selling reactive power in the following market. This makes it difficult to leverage DG’s ability to operate at low PFs as they will be incentivized to sell real power first, thus limiting reactive power ability.

The third category covers simultaneous real and reactive power deployment. Typical examples are \cite{yang_linearized_2018,bai_distribution_2018,el-samahy_procurement_2008} that extract both real and reactive power prices using an OPF-formulation. By minimizing an objective function, the resulting Lagrange multipliers (dual variables) represent the marginal cost of real and reactive power at each load and resource, leading to a reactive power price. For example, \cite{yang_linearized_2018} and \cite{bai_distribution_2018} use voltage and reactive power constraints to ensure grid-feasible solutions and find the value of reactive power in meeting these constraints. Several studies have built on this approach to address more practical concerns such as  \cite{roozbehani_dynamic_2010} which proposes an algorithm to limit and characterize price volatility and \cite{moghimi_short-term_2020} which divides the market at the transmission level into profit-seeking and a social-welfare seeking components. By solving for both types of power simultaneously, these approaches can avoid decoupling grid dynamics and assure there is no mismatch between cleared real and reactive dispatch. Our proposed market solution belongs to this third category.



\subsection{Reactive power needs of the emerging grid}
The state-of-art in reactive power management and pricing are insufficient in accommodating a deep penetration of DERs, which is likely to occur in the future grid. In order to leverage the flexibility and low cost of DERs, a new set of requirements and objectives are needed. For the emerging low-carbon grid with high DER penetration, we propose the following set of new priorities for industry practice:
\begin{itemize}
    \item Utilize the significant reactive power capabilities of a large penetration of DERs and their wider PF operating range 
    \item Support efficient grid operation by minimizing grid losses using local DERs
    \item Provide fair remuneration to DERs for their grid services
    \item Incentivize long-term growth of DER adoption through market inclusion and corresponding reliable revenue streams
\end{itemize}
In the next section, we introduce our reactive power market and pricing mechanism for the emerging grid.


\section{Reactive Power Pricing} \label{sec:market-design}

\begin{figure}[h!]
    \centering
    \includegraphics[scale=0.5]{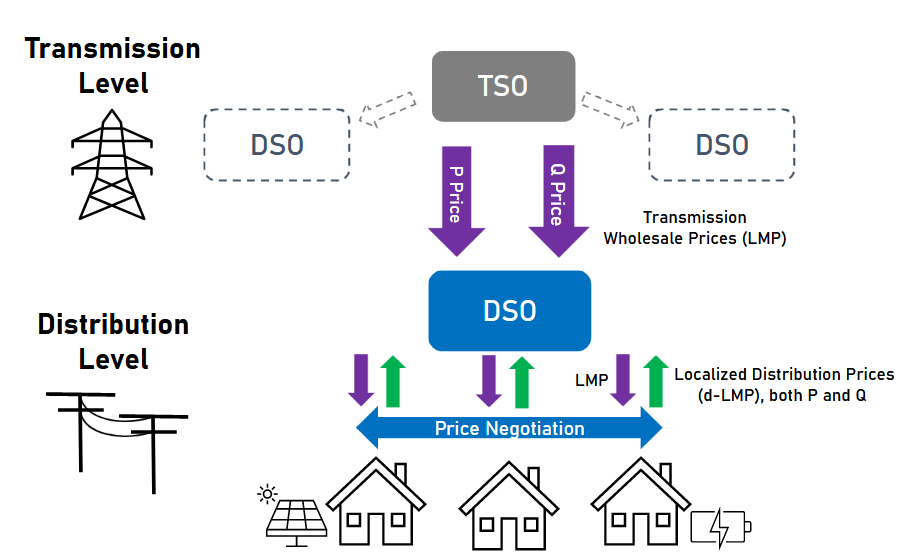}
    \caption{The proposed retail market structure is outlined here. The DSO oversees the market operation. Using wholesale LMPs at the distribution grid feeder and other grid information, the market simultaneously clears locational prices for real and reactive power. Price negotiations occur between resource owners using the PAC algorithm.}
    \label{fig:market_structure}
\end{figure}

The starting point for the reactive power market is an unbalanced distribution grid for which we are interested in determining the optimal power dispatch, for both active and reactive power simultaneously. Figure~\ref{fig:market_structure} illustrates how the TSO and DSO interact; using the wholesale prices for real and reactive power at the distribution substation, the DSO facilitates price negotiations to determine local and time-varying prices for real and reactive power. We note that while reactive prices at the wholesale level do not currently exist, we include a Q-LMP for completion. We utilize an OPF approach to solve this problem, whose primal and dual variables provide the power dispatch and the local price of both active and reactive power. These local prices vary temporally and spatially, providing more accurate economic signals for DER control and fair service-based remuneration. It is worth noting this new energy market provides tertiary control through price incentives for the grid, and may be viewed as a transactive energy scheme\footnote{Defined by NIST as ``a system of economic and control mechanisms that allows the dynamic balance of supply and demand across the entire electrical infrastructure using value as a key operational parameter'' \citep{NISTtransactive}, transactive energy bridges the gap between physical power flow in the grid and market derivatives.}. The proposed market does not replace low-level voltage controllers tasked with primary voltage control, just as real power energy markets do not replace frequency control. The design of ancillary markets and reactive power pricing for voltage control is beyond the scope of this paper.

Within the distribution grid, many DERs may be connected to only a single phase of the three-phase network, necessitating a high-fidelity power flow model of the unbalanced distribution grid. In what follows, we present a detailed discussion of the current injection model of the unbalanced grid, a distributed optimization based market mechanism, and the specifics of reactive power pricing.

\subsection{Introduction to Current Injection Model}
Two distinct types of models have been used to represent the physics of the distribution grid: modeling branch variables which leads to the Branch Flow model \cite{lavaei_2012,farivar_2014}, or modeling nodal variables which leads to the Bus Injection model \cite{Low2014_p1,emiliano_2012}. The Branch Flow model based on second order cone programming (SOCP) has proven to be advantageous in providing tight convex relaxations to the original AC-OPF problem with exactness under some conditions, and is shown to be more computationally stable than the Bus Injection model. However, both Branch Flow and Bus Injection models are typically limited to networks with radial topologies and balanced networks, and extensions to these models for unbalanced distribution grids are only valid for a small range of angle imbalances. This is a limiting assumption for many grids, especially with increasing penetration of DERs located on single-phase lines. A recent approach proposed by some of the authors of this paper, denoted as Current Injection (CI) model \cite{FerroThesis,FerroCIM_IFAC}, avoids this assumption and so is an ideal candidate for representing unbalanced grids with various single-phase loads and generation. The CI model uses nodal variables, similar to the Bus Injection model, but the main idea is to represent all loads and generators as nodal current injections, with all power, current, and voltage phasors represented in Cartesian coordinates. The 3-phase impedance matrix is used to describe the self and mutual inductance between phases to model the coupling of phases that are common to a distribution grid. More importantly, the key obstacle of non-convexity of the AC-OPF and the subsequent nonlinearity of SOCP and SDP convexification strategies are dealt with in the CI approach by leveraging McCormick Envelope (MCE) based convex relaxation \cite{mccormick_1976} for the bilinear power relations. The MCE uses the convex hull representation of bilinear terms to render a linear OPF model. This representation requires adequate bounds on the nodal voltages and currents to ensure a tight convex relaxation. To determine these bounds, the CI approach also includes a carefully designed pre-processing algorithm which uses generation and load forecasts and grid limits to iteratively calculate tight nodal bounds \cite{FerroThesis}. The CI model has been shown to perform well on unbalanced networks with local generation, with maximum 1.2\% optimality gap and 0.9\% voltage error when compared to the AC-OPF for a number of use cases \cite{FerroThesis}. The CI model has also been used for solving OPF in a distributed way \cite{FerroCIM_IFAC} and for voltage support from inverter-based resources in unbalanced distribution grids \cite{Haider_VVC}. Due to its overall ability to model unbalanced grids and all of the aforementioned advantages including the computational simplicity of the linear model, we adopt the CI approach in this paper. 

\subsection{Problem Formulation} \label{sec:market-design:CI}
We denote a general distribution network as a graph $\Gamma(\mathcal{N},\mathcal{E})$, where $\mathcal{N}\coloneqq\{1,...,N\}$ denotes the set of nodes, $\mathcal{E}\coloneqq\{(m,n)\}$ denotes the set of edges, and each phase is expressed as $\phi\in \mathcal{P}, \; \mathcal{P}=\{a,b,c\}$. The CI-OPF for the network is then written as:
\begingroup
\allowdisplaybreaks
\begin{subequations} \label{eq:CI}
\begin{align}
    \underset{x = \left[{{I}^{R}}\;{{I}^{I}}\;{{V}^{R}}\;{{V}^{I}}\;P\;Q\right]}{\mathop{\min }}\; f & \left(x\right) \label{eq:CIOPF_obj} \\
    V &= Z I \label{eq:CIOPF_ohm} \\
    P_{j}^{\phi} &= V_{j}^{\phi,R}I_{j}^{\phi,R}+V_{j}^{\phi,I}I_{j}^{\phi,I}  \label{eq:CIOPF_Pdef} \\
    Q_{j}^{\phi} &= -V_{j}^{\phi,R}I_{j}^{\phi,I}+V_{j}^{\phi,I}I_{j}^{\phi,R}  \label{eq:CIOPF_Qdef} \\
    P_{j}^{\phi} \tan(\cos^{-1}(-\text{pf})) &\leq Q_j^{\phi} \leq P_j^{\phi} \tan(\cos^{-1}(\text{pf})) & \label{eq:CIOPF_pf} \\
    \underline{P}_{j}^{\phi} &\le P_{j}^{\phi} \le \overline{P}_{j}^{\phi}  \label{eq:CIOPF_Plim} \\
    \underline{Q}_{j}^{\phi} &\le Q_{j}^{\phi}\le \overline{Q}_{j}^{\phi} \label{eq:CIOPF_Qlim} \\
    \underline{V}_{j}^{\phi,R} &\le V_{j}^{\phi,R} \le \overline{V}_{j}^{\phi,R}  \label{eq:CIOPF_Vlim1}\\
    \underline{V}{j}^{\phi,I} &\le V_{j}^{\phi,I} \le \overline{V}_{j}^{\phi,I} \label{eq:CIOPF_Vlim2} \\
    \underline{I}_j^{\phi,R} &\le I_j^{\phi,R}\le \overline{I}_j^{\phi,R} \label{eq:CIOPF_Ilim1} \\
     \underline{I}_j^{\phi,I} &\le I_j^{\phi,I}\le \overline{I}_j^{\phi,I}  \label{eq:CIOPF_Ilim2}
\end{align}
\end{subequations}
\endgroup
where $x = \left[{{I}^{R}}\;{{I}^{I}}\;{{V}^{R}}\;{{V}^{I}}\;P\;Q\right]$ is the decision vector for the CI-OPF problem; $I, V, P, Q$ denote the vector of nodal current injections, voltages, and real/reactive power injections respectively; $Z\in \mathbb{R}^{3N\times3N}$ is the system impedance matrix, which includes cross-phase coupling inherent in multi-phase distribution grids. We use $x^R$ and $x^I$ to denote the real and imaginary components of a complex number $x$; overbar $\overline{{x}}$ and underbar $\underline{{x}}$ denote the upper and lower limits of a variable ${x}$; $\operatorname{Re}(\cdot)$ and $\operatorname{Im}(\cdot)$ denote the real and imaginary components of a complex number. Constraint \eqref{eq:CIOPF_ohm} describes Ohm's law, and \eqref{eq:CIOPF_Pdef}-\eqref{eq:CIOPF_Qdef} are the definitions of real and reactive power. Constraint \eqref{eq:CIOPF_pf} models the multiphase inverter of a DG, where the ratio of P and Q determine the PF setting. Constraints \eqref{eq:CIOPF_Pdef}-\eqref{eq:CIOPF_Ilim2} are for all nodes $j \in \mathcal{N}$ and per each phase $\phi \in \mathcal{P}$. 

\textit{Remark 1:} The use of current and voltage to model the power physics permits a simple multi-phase model for the unbalanced distribution grid. The corresponding definition of real/reactive power as represented in constraint \eqref{eq:CIOPF_Pdef}-\eqref{eq:CIOPF_Qdef} include bilinear terms which can be nonconvex. To deal with this nonconvexity, we use a McCormick Envelope based convex relaxation, which renders a linear model. We then convert the bounds on P and Q to nodal current limits, $\underline{I}_j^{\phi}$ and $\overline{I}_j^{\phi}$, for both the real and reactive components, using a carefully designed pre-processing algorithm. For brevity, we omit the McCormick Envelope relaxation of the bilinear terms and the pre-processing. See \cite{FerroCIM_IFAC,FerroThesis} for details.

\textit{Remark 2:} The CI model includes hard constraints on the nodal voltages in \eqref{eq:CIOPF_Vlim1} and \eqref{eq:CIOPF_Vlim2}. These bounds on the real and imaginary components of the voltage phasor at every node are determined by the pre-processing technique, using information of operating bounds (ex. $\pm 5\%$ in the US) and the forecast for power injections and loads. Thus, any dispatch solution to the CI-OPF model will enforce voltage constraints, as required for grid operation.

\textit{Remark 3:} Another typical constraint in power systems is a congestion constraint, which limits the line currents. This constraint can be trivially included in our CI-OPF formulation by introducing a variable describing the current flow in the distribution lines, $I_\text{flow}$, and rewriting constraint \eqref{eq:CIOPF_ohm} as two constraints: 
\begin{subequations}
\begin{align} \label{eq:Iflow}
    AV &= ZI_{\textrm{flow}} \\
    I &= A^TI_{\textrm{flow}}
\end{align}
\end{subequations}
where $A \in \mathbb{R}^{3N\times3N}$ is the 3-phase graph incidence matrix.

The objective function $f(x)$ in \eqref{eq:CIOPF_obj} is chosen to carefully reflect the energy market structure. As previously discussed, our proposed market solution covers simultaneous real and reactive power deployment, and energy price settlements. To reflect these energy prices, the objective function is chosen to be a combination of Social Welfare and line losses and is given by:
\begin{align}
    f \left(x\right) = &\sum_{j\in \mathcal{N}_G} \left( a_{j}^{P}P_{j}^2 + b_{j}^{P}P_{j} + a_{j}^{Q}Q_{j}^2 + b_{j}^{Q}Q_{j} \right) \nonumber \\
    & - \sum_{j\in \mathcal{N}_L} \left(\alpha_{j}^{P} (P_{j}-\underline{P_{j}})^2
    + \alpha_{j}^{Q}(Q_{j}-\underline{Q_{j}})^2\right) \nonumber \\
    & + \zeta \operatorname{Re}\left(I^{H}ZI\right) \label{eq:CIOPF_obj_econ} 
\end{align}
where the first two terms represents generator cost and customer load disutility, and the third term represents line losses, in which superscript $H$ denotes the Hermitian. In \eqref{eq:CIOPF_obj_econ}, $\mathcal{N}_G$ and $\mathcal{N}_L$ are the set of generator and load nodes respectively, $\zeta$ describes the tradeoff between economic and energy efficiency objectives, $a_j^P, b_j^P,a_j^Q,b_j^Q$ are generating cost coefficients, and $\alpha_j^P, \alpha_j^Q$ are load disutility coefficients. All nodes $j\in \mathcal{N}_G$ and $j\in \mathcal{N}_L$ represent active DG and DR participants in the distribution grid, respectively. Node 1 is treated as the point of common coupling with the transmission grid, and reflects the wholesale price of electricity, which is the LMP $\lambda^P$ from the WEM. Thus for $j = 1$, $a_j^P = 0$ and $b_j^P = \lambda^P$. The cost coefficients related to reactive power are chosen as $a_j^Q=0,\quad b_j^Q=0.1 b_j^P$. The motivation for these choices comes from \cite{federal_energy_regulatory_commission_payment_2014} which cites that reactive power prices are often one-tenth of that of real power. It should be noted that the weighted combination of the social welfare and line losses is not standard practice, but is included here as line losses are more significant for an optimal functioning of a distribution grid compared to a transmission grid. Further analysis of this objective function as it relates to reactive power pricing is presented in Section~\ref{sec:market-design:pricing}. The resulting solutions of the CI-OPF problem forms the backbone of our proposed reactive power market.


\subsection{Distributed Market Mechanism} \label{sec:market-design:algorithm}
Our proposed retail market mechanism enables market participation of a large number of DERs spatially distributed throughout the grid, and involves simultaneous P-Q price and device set point calculations. To solve the corresponding CI-OPF model in a centralized manner can become quickly intractable with a large number of DERs participating. An attractive alternative to centralized perspectives is a distributed one. Herein, each agent has access to local information and a limited amount of information shared by its neighbours, which it uses to arrive at the market solution through iterative computation and communication. Such a distributed approach leverages the computational and communication capabilities of intelligent grid edge devices. Further, if well designed, they can be resilient to communication link and single-point failures, and can preserve the private information of the DERs while still realizing network-level objectives. Distributed decision making is an emerging paradigm in the distribution grid and can enable the decarbonization of future energy systems. Thus, we take a distributed optimization approach for our retail market. In what follows, we introduce a Proximal Atomic Coordination (PAC) based algorithm for our retail market.

The CI-OPF problem \eqref{eq:CI} can be cast into standard form as:
\begin{gather}
    \underset{x}{\mathop{\min }} \;\;  \frac{1}{2}x^T M x + c^Tx \nonumber\\
    \text{s. t.} \begin{cases}
        Gx=b \\
        Hx\le d
    \end{cases}
    \label{eq:general_problem}
\end{gather}
where $M$ and $c$ are the quadratic and linear cost matrix and vector respectively, $G$ and $b$ are the equality constraint matrix and vector (i.e. representing \eqref{eq:CIOPF_ohm}-\eqref{eq:CIOPF_Qdef}), and $H$ and $d$ are the inequality constraint matrix and vector (i.e. representing \eqref{eq:CIOPF_pf}-\eqref{eq:CIOPF_Ilim2}).

\subsubsection{Decomposition of CI-OPF}
Problem \eqref{eq:general_problem} can then be reformulated into $N$ separate coupled optimization problems, each of the form:
\begin{gather}
    \underset{a}{\mathop{\min }} \;\; \sum_{j\in \mathcal{N}} \frac{1}{2}a_j^T M_j a_j + c_j^Ta_j \nonumber\\
    \text{s. t.} \begin{cases}
        G_j a_j = b_j  & (\mu_j)\\
		H_j a_j \le d_j  & (\lambda_j)\\
		B_j a = 0  & (\nu_j)
    \end{cases}
    \label{eq:general_atomized}
\end{gather}
where $M_j, c_j, G_j, b_j, H_j, d_j$ represent the local cost and constraint submatrices and subvectors, and $\mu_j, \lambda_j, \nu_j$ represent the corresponding dual variables for each constraint. In order to fully distribute the coupled CI-OPF, each agent $j$ will \textit{estimate} the value of coupling variables which are \textit{owned} by a neighbouring agent $i$. A set of \textit{coordination constraints} must be added for each subproblem $j\in N$:
\begin{equation}
	 B_j {a} = 0 \label{eq:coord_const}
\end{equation}
where $B$ represents the incidence matrix of a directed graph whose nodes are the owned and copied atomic variables within $a_{1}, \ldots, a_{K}$.  Hence, we have for each element of the incidence matrix (where $i$ and $m$ correspond to indices of atomic variables):
\begin{equation*}
    B_{i}^{m} \triangleq \left\{\begin{array}{rl} -1, & \text{if $i$ is `owned' and $m$ the related `estimate'} \\ 
    1, & \text{if $m$ is `owned' and $i$ a related `estimate'} \\
    0, & \text{otherwise}\end{array}\right.
\end{equation*}
Submatrices ${B_j}$ and $B^j$ represent the relevant incoming and out-going edges of the directed graph for subproblem $j$, respectively. This coordination constraint \eqref{eq:coord_const} can be interpreted as requiring all subproblem estimated variables in a given ${j}$th subproblem  to equal the value of their corresponding owned in $i$th subproblem, $i\neq j$.
In the CI-OPF, the coupling comes from the power physics of the grid, namely constraint \eqref{eq:CIOPF_ohm}. The resulting atomic decision variable for agent $j$ is:
\begin{equation}
    a_j = \left[ I_j^R ; I_j^I ; V_j^R ; V_j^I ; P_j ; Q_j ; \{\hat{I}_{i}^R,\hat{I}_{i}^I\}_{\forall i, \{(i,j)\} \in \mathcal{E}} \right] \label{eq:atomized_var}
\end{equation}
where variable estimates are denoted with a hat.

\subsubsection{Statement of the Algorithm}
\noindent The Lagrangian function for each network subproblem is:
\begin{align}
    \mathcal{L}_{{j}} \left({a}_{{j}}, {\mu}_{{j}}, {\nu}, \lambda_j \right) = \sum_{{j} \in K} &\left[{f}_{{j}} \left({a}_{{j}}\right) + {\mu}_{{j}}^{T} ({{G}}_{{j}} {a}_{{j}}-b_j) \right. \nonumber\\
    & \left.+ {\nu}^{T}  B^j {a}_{j} +\lambda_{{j}}^{T}\left ({{H}}_{{j}} {a}_{{j}}  -{d}_{j}\right  )\right] \label{eq:lagrange}
\end{align}
where it can be shown that 
\begin{equation}
    \mathcal{L} \left({a},{\mu},{\nu}, \lambda \right) \triangleq \sum_{{j} \in K} \mathcal{L}_{{j}} \left({a}_{{j}}, {\mu}_{{j}}, {\nu}, \lambda_j \right) \nonumber 
\end{equation}

The PAC-based distributed algorithm to solve the CI-OPF is then given by:
\begin{subequations} \label{eq:PAC}
\begin{align}
    {a}_{{j}} \left[ \tau + 1 \right] &= -\left[M_j +\rho\gamma \left(G_{j}^{T}G_{j}+B_{j}^{T}B_{j}+\frac{1}{\rho}\right)\right]^{-1} \nonumber \\
    ( c_j + & {\mu}_j^T  G_j + {\nu}^{T}  B^j +\lambda_j^T H_j -\rho \gamma G_j b_j -\frac{a_j \left[ \tau \right] }{\rho} ) \\
    \lambda_j \left[ \tau + 1 \right] &= \mathop{\max} \left[ 0,\, \lambda_j \left[ \tau \right] + \rho \gamma \left ( {{H}}_{{j}} {a}_{{j}} \left[ \tau+ 1\right]-d_j\right ) \right] \\
    \mu_j \left[ \tau + 1 \right] &= \mu_j \left[ \tau \right] + \rho \gamma \left ( G_j a_j \left[ \tau+ 1\right]-b_j\right ) \\
	\bar{\mu}_j \left[ \tau + 1 \right] &= \mu_j \left[\tau + 1\right] +\rho \gamma_j[\tau] \left ( G_j a_j \left[ \tau+ 1\right]-b_j\right ) \\
	\text{Commun} & \text{icate } \left\{a_j\right\}_{j \in N} \text{ within network} \\
	\nu_j \left[ \tau + 1 \right] &= \nu_j \left[ \tau \right] + \rho \gamma B_j a\left[\tau+1\right]\\
	\bar{\nu}_j \left[ \tau + 1 \right] &= \nu_j \left[ \tau + 1\right] + \rho \gamma_j[\tau] {B}_j {a}\left[\tau+1\right] \\
	\text{Commun} & \text{icate } \left\{\bar{\nu}_j\right\}_{j \in N} \text{ within network}
\end{align}
\end{subequations}
The PAC algorithm achieves convergence in both objective value and distance to feasibility with a linear rate for convex problems - such as that of the CI-OPF. 

\begin{theorem}
    Let ${a} \left[\tau\right] = \left[{a}_{{1}}\left[\tau\right]; \cdots; {a}_{{K}}\left[\tau\right]\right]$
	represent the PAC trajectory of \eqref{eq:PAC}. Define the following quantities:
\begin{equation*}
    	\tilde{{G}} = \text{{diag}} \left\{\tilde{{G}}_{{1}}, \ldots, \tilde{{G}}_{{K}}\right\}
\end{equation*}
\begin{equation*}
    {V}_{{1}} = \tilde{{G}}^{T} \tilde{{G}} + {B}^{T}{B},
\end{equation*}
\begin{equation*}
    	\tilde{{V}}_{{1}} \left({\gamma}\right) =\gamma{V}_{{1}}
\end{equation*}
\begin{equation*}
    \tilde{{V}}_{{2}} \left(\rho, \gamma\right) = \frac{1}{\rho^{2}} {I}_{{N}^{\mathsf{T}}} - \tilde{{V}}_{{1}} \left({\gamma}\right)
\end{equation*}Then if the following assumptions hold
\begin{itemize}
    \item The objective function $f(x)$ is convex;
    \item The PAC parameters satisfy		$1 > \rho^{2} \gamma \lambda_{\text{{max}}} \left(\tilde{{G}}^{T} \tilde{{G}} + {B}^{T} {B}\right)$ with  $\rho > 0$ and $\gamma> 0$.
\end{itemize}
Then there exists an optimal solution ${a}^{\ast}$ s.t.:
		\begin{equation*}
			\lim_{\tau \rightarrow \infty} \left\{{a} \left[\tau\right]\right\} = {a}^{\ast},
		\end{equation*}
		with convergence rate satisfying, for all $\tau \in \mathbb{N}$:
		\begin{equation*}
			\left\{\left\Vert {a} \left[\tau + 1\right] - {a} \left[\tau\right] \right\Vert_{{V}_{{2}}\left(\rho,{\gamma}\right)}\right\}_{\tau \in \mathbb{N}} = o\left(\frac{1}{\tau}\right).
		\end{equation*}
\end{theorem}

The proof of convergence results and advantages with respect to state-of-the-art solutions for the PAC-based distributed algorithm are detailed in \cite{RomvaryTAC}.

\subsubsection{Privacy of the Algorithm}
The distributed PAC algorithm has enhanced privacy preserving characteristics. We consider sensitive information to correspond to any information about the ${k}$-th market participant that the ${j}$-th participant can make use of to their advantage -- such as cost functions and operating constraints. This protection of this sensitive information corresponds to privacy. Within the PAC algorithm, the objective function, constraints, and dual variables $\mu_j$ are not shared between the $j$-th and $k$-th atom. Further, the dual variables $\nu_j$ are not shared directly; rather, a ``protected'' Lagrange multiplier $\bar{\nu}_j$ is shared with neighbours. The knowledge of this dual variable can be used to sabotage the convergence of PAC (and the corresponding market clearing) through deliberate manipulations of the dual variables corresponding to the coordination constraints. By using $\bar{\nu}_j$, the PAC algorithm provides a privacy guarantee to protect against this form of data manipulation, and ensures the market can clear.

Denoting the rogue atom's estimation is carried out as:
\begin{equation}
    \hat{\nu }_{j\to k}\left[ \tau  \right] = \rho \widehat{{\gamma}_{j} B_j} \left[\sum_{s = 0}^{\tau} {a} [s]\right] , \label{eq:estimation}
\end{equation}
where $\widehat{\gamma_{j}B_j}$ is an estimate of ${{\gamma}_{{j}} B_j}$. By using a time-varying over-relaxation rate $\hat{\gamma}_j\left[\tau\right]$, the underlying ${\gamma}_{j} B_j$ matrix cannot be recovered, as given by Proposition \ref{prop_privacy} from \citep{RomvaryTAC}, and repeated below:
\begin{proposition} \label{prop_privacy}
Given the algorithm \eqref{eq:PAC}, when the rogue agent estimation is carried out as in \eqref{eq:estimation}, then
\begin{equation}
    \| \hat{\nu }_{{j}\to {k}}\left[\tau \right] -{\nu}_{j} [\tau] \|_2 \not\rightarrow 0 \label{eq:proof_privacy}
\end{equation}
\end{proposition}

\subsection{Reactive Power Pricing} \label{sec:market-design:pricing}
The prices for real and reactive power of a DER at node $j$ are chosen as the marginal costs of generators/consumers and are therefore determined using the dual variables of the CI-OPF model. In particular, the Lagrange multipliers corresponding to \eqref{eq:CIOPF_Pdef} and \eqref{eq:CIOPF_Qdef} are used to determine the price for real and reactive power respectively. While the discussions below are focused almost entirely on the reactive power component, the solutions above encompass both the real and reactive power dispatch simultaneously.

Each market clearing sets the reactive power dispatch, $Q_{j,t}$ for node $j$ and time $t$, and produces the corresponding dual variable $\mu_{j,t}^Q$. This dual variable sets the basis for the time-varying reactive power price $\bar{\mu}_j^Q$, which we denote as the d-LMP, determined for node $j$ as:
\begin{align}
    \bar{\mu}_j^Q = \frac{\sum_{t \in T} \mu_{j,t}^Q Q_{j,t}}{\sum_{t \in T} Q_{j,t}}
    \label{eq:price_def}
\end{align}
This d-LMP varies daily, calculated as a weighted average of the dual variable from the OPF problem in \eqref{eq:CI}. Thus for a 1-hour retail market clearing period, the set $T$ includes 24 points. A more frequent clearing can also be accommodated by our proposed market structure, particularly when forecast errors are substantial. Such a time-varying price allows DGs to adjust their generation and/or PF settings and DRs to shift their consumption behavior, all in a coordinated manner so that the DSO can accurately recover costs. Rather than exposing end-use customers and DER owners to the complete dynamics of the electricity system by using $\mu_{j,t}^Q$ as a real-time d-LMP, as wholesale market participation models (like FERC 2222) would do with a corresponding $\lambda^Q$, the averaging procedure in \eqref{eq:price_def} allows the price volatility to be contained while resulting in the same payout as the corresponding real-time prices.

\subsubsection{Reactive price volatility}
An important and necessary distinction here is made between the dispatch decisions and financial transactions within the reactive power market. Suppose each DER participates in a 5-minute retail market. For each market period, DERs will submit iterative bids to neighbouring agents (as per the distributed PAC-based algorithm), then upon agreement (i.e. convergence of the algorithm), clear the market at the agreed-upon real-time price $\mu_{j,t}^Q$ and commit to dispatch $Q_{j,t}$. This allows the market to use 5-minute ahead forecasts for DERs, thereby significantly reducing forecast error as compared to a 24-hour ahead window of a day-ahead market, while creating the requisite signal for DERs to respond to grid needs at a faster timescale than an hourly or daily signal. Meanwhile, the financial transaction, i.e. a customer's bill or payout to the DER, will reflect the daily d-LMP $\overline{\mu}_j^Q$ as calculated in \eqref{eq:price_def}. In using a weighted average, the perception of price volatility is removed, instilling confidence in the retail market. This simple approach to volatility reduction is shown to be sufficient through simulations in Section~\ref{sec:results}, where our reactive power market results in consistent revenue for DER owners over the week.

\subsubsection{Market objective function}
In this section we present an analysis of how the objective function \eqref{eq:CIOPF_obj_econ} affects the individual reactive power prices at the nodes. The first two terms of the objective function represent the generator cost and customer load disutility, capturing the P and Q prices at a subset of nodes, i.e. generators and flexible loads. The third term represents line losses (${\operatorname{Re}\left(I^{H}ZI\right)}$) and captures a global view of the network, introducing Q prices at every node. In what follows, we show this relation between line losses and Q prices.

The losses in an AC network are calculated as:
\begin{align}
   S_l &= \left( I\right)^H Z I \nonumber\\
   &=\left( I^{I}\right)^T R I^{I} +\left( I^{R}\right)^T R I^{R} \nonumber\\ & +  i\left[\left( I^{I}\right)^T X I^{I} + \left( I^{R}\right)^T X I^{R}\right]
\end{align}
where $i$ is the imaginary unit. The corresponding real and imaginary part of the complex number $S_l$ define the real and reactive power losses, with the objective function \eqref{eq:CIOPF_obj_econ} minimizing the real power loss as the third term:
\begin{equation}
   P_l= \left( I^{I}\right)^T R I^{I} +\left( I^{R}\right)^T R I^{R} \label{eq:P_loss}
\end{equation}

Recall that the nodal prices are defined as the Lagrange multipliers associated with Eqs. \eqref{eq:CIOPF_Pdef}-\eqref{eq:CIOPF_Qdef}. These constraints can be re-written in terms of the nodal currents as:
\begin{align}
    Q = & -V^RI^I + V^II^R \nonumber \\
    = & -(RI^R-XI^I)I^{I}+(RI^I+XI^R)I^{R} \nonumber \\
    = & -I^{I} \odot RI^R+I^{I}\odot XI^I+I^{R}\odot RI^I+I^{R}\odot XI^R \nonumber \\
    = & I^{I}\odot XI^I+I^{R}\odot XI^R \label{eq:Q_fromcurrents}
\end{align}
where $\odot$ is the Hadamard product, and using the relationship $V^R=RI^R-XI^I$ and $V^I=RI^I+XI^R$ from \eqref{eq:CIOPF_ohm}. Similarly, we can show $P = I^{R} \odot RI^R+I^{I}\odot RI^I$. Comparing \eqref{eq:P_loss} and \eqref{eq:Q_fromcurrents}, the reactive power $Q$ is a nonlinear combination of the variables composing the loss term $P_l$. By minimizing $P_l$, the objective function then directly influences the nodal prices of reactive power. Thus, minimizing power loss serves as a network-wide view of reactive power, and generates an appropriate price signal. It is important to note that by virtue of Eq.~\eqref{eq:CIOPF_ohm}, the minimization of losses in the objective function implicitly minimizes the voltage drop across the whole network.

In the next section we present a numerical case study of the reactive power market and pricing model.

\section{Method}  \label{sec:method}

Our proposed reactive power market for the distribution grid is evaluated in this section using an IEEE 123-bus network. The network is modified to include clusters of rooftop photovoltaic (PV) units at up to 27 of the 123 nodes. To emulate varying levels of DG penetration, PV clusters are incrementally added at various nodes. Each PV cluster is chosen to have a capacity of 25-80kW, a range that represents an addition of 3 to 10 small residential installations or larger projects such as carports or community solar. 

In our simulations, the cluster capacity range was chosen so as to emulate a range of 5\% to 160\% DG penetration (ratio of nameplate capacity to average network load, see \eqref{eq:penetration}), which can reflect the trend that is anticipated over the coming decade especially in the context of more aggressive DG adoption scenarios. For New England (US), the current and forecasted DG penetration by 2030 is 4.3\% and 14\%, respectively \cite{iso_ne_concentration}. Each DG is assumed to be equipped with a smart inverter which operates at a flexible PF, with appropriate power electronic control \cite{ieee_DER_interconnection}.
Rather than fixing the PF, we enforce a lower limit for each unit, permitting the market to determine the inverter setting. We model all inverters as having the same PF limits. The baseline PV generation data $\alpha^\text{PV}(t)$ is taken from the NREL SAM dataset \cite{nrel_dataset}, and is assumed to be the same for each PV cluster. The load profile for each node is modelled with real-time data from ISO-NE \cite{iso_express}, and perturbed for each load at node $j$ as $\alpha_j^P(t)=\alpha_D(t)\delta^P_{j}$ with $\delta^P_{j} \sim \mathcal{N}(0,0.1)$. Wholesale electricity prices are obtained from ISO-NE's real-time market to determine $\lambda^P$. We assume our market clears hourly. With these modifications introduced into the IEEE-123 bus network, the reactive power price \eqref{eq:price_def} and the reactive power injection $Q_j$ are determined using the PAC-based distributed CI-OPF described by Eq.~\eqref{eq:PAC}.



We assess the performance of our proposed retail market by considering the following scenarios:
\begin{enumerate}
    \item with increasing PV penetration; and
    \item with varying lower limit on PF operating range
\end{enumerate}

We also quantify the performance of our retail market by evaluating the variation in reactive power pricing over a 24-hour period, and over a week. To assess the performance of our reactive power market, we introduce the following metrics. 

\subsection{Scenario metrics}
We define DG penetration as 
\begin{align}
    \text{DG penetration} = \frac{\sum_{ \mathcal{N}_G} \textrm{Nameplate Capacity}}{\sum_{\mathcal{N}_L} \textrm{Average Load}} \label{eq:penetration} 
\end{align}
We also introduce a second penetration metric which accounts for the capacity factor of renewable resources, the DG energy penetration, as
\begin{align}
    \text{DG energy penetration} = \frac{\sum_{t \in T}\sum_{j \in \mathcal{N}_G}S_{j,t}}{\sum_{t \in T}\sum_{j \in \mathcal{N}_L}S_{j,t}} \label{eq:en_penetration} 
\end{align}
where $S$ is apparent power. We further note that as the capacity factors of DGs are often small, the \textit{DG penetration} metric may exceed 100\% to reach a higher \textit{DG energy penetration}, ensuring carbon-free solar energy is available to serve load on the feeder. This corresponds to the installed nameplate capacity exceeding the average network load on the feeder. Further, community solar projects may have a larger capacity than demand on the same feeder, as they may have shares owned by customers connected to adjacent feeders, or by commercial/private owners investing in solar projects.

\subsection{Revenue metrics}
\noindent The ratio of the reactive power revenue to the total revenue for DGs is defined as
\begin{align}
    \text{Q-} & \text{Revenue Ratio} =  \frac{1}{|\mathcal{N}_G|} \sum_{j \in \mathcal{N}_G} \frac{\left(\bar{\mu}_j^Q\sum_{t\in T} Q_{j,t} \right)}
    {\left(\bar{\mu}_j^Q\sum_{t\in T} Q_{j,t} + \sum_{t \in T}\mu^{P}_{j,t}P_{j,t}\right)}
 \label{eq:revenue_ratio}
\end{align}
where $|\cdot|$ denotes the cardinality of set $\mathcal{N}_G$, i.e. the number of generator nodes. That is, the Q-revenue ratio is the average fraction of the reactive power payout compared to the total payout from real and reactive power, over the network. 

\noindent The coefficient of variation $c_v$ is defined below as
\begin{align}
    c_v = \frac{std \left( \mu_j^Q \right)}{ \overline{\mu_j^Q}}
    \label{eq:coeff_variation} 
\end{align}
and is a measure of the reactive price volatility.

\subsection{Grid performance metrics}
\noindent The DG-Q utilization is defined as
\begin{align}
     \text{DG-Q utilization} = \frac{\sum_{t \in T}\sum_{j \in \mathcal{N}_G}Q_{j,t}}{\sum_{t \in T}\sum_{j \in \mathcal{N}_L}Q_{j,t}}
  \label{eq:Qpenetration} 
\end{align}
That is, the DG-Q utilization is the fraction of reactive power loads that are served by DGs. Similarly, the DG-P utilization is defined as 
\begin{align}
     \text{DG-P utilization} = \frac{\sum_{t \in T}\sum_{j \in \mathcal{N}_G}P_{j,t}}{\sum_{t \in T}\sum_{j \in \mathcal{N}_L}P_{j,t}}
  \label{eq:Ppenetration} 
\end{align}
The remaining load (for both P and Q) is assumed to be served by the transmission system. 

\noindent The network line losses are defined as the sum of losses over the grid calculated as
\begin{align}
    \text{Network Losses} =  \sum_{j \in \mathcal{N}_L} R_{j} \left( I_{j}^{I^{2}}+I_{j}^{R^{2}} \right)
 \label{eq:network_losses}
\end{align}

\noindent where $R_j$ is the resistance of line $j$. The Mean Voltage is the voltage averaged over nodes and time defined as
\begin{align}
    \text{Mean Voltage} =  \frac{1}{|\mathcal{N}_G| \times T}\sum_{j \in \mathcal{N}_G}\sum_{t\in T} V_{j,t}
 \label{eq:mean_voltage}
\end{align}

\noindent The PF per DG $k$ is calculated as
\begin{equation*}
    PF_k = \cos\left(\arctan\frac{\sum_{\phi \in \mathcal{P}}Q_k^{\phi}}{\sum_{\phi \in \mathcal{P}}P_k^{\phi}}\right) \label{eq:pf}
\end{equation*}

\section{Results} \label{sec:results}
\subsection{Market Response to Increasing DG Penetration}
In this section we investigate the behaviour of the reactive power market for a range of DG penetration of 5\% to 160\% obtained from the modified IEEE-123 bus. This corresponds to a DG energy penetration range of 0.5\% to 35\%, meaning 35\% of the total daily load is served by local PV generation and the remaining serviced by imports from transmission grid. All inverters were assumed to have a PF limit of 0.9. By monitoring Q prices and production as DGs are added, we show that our market is capable of utilizing and incentivizing Q-generation from DGs at a high penetration.

\begin{figure}[h!]
    \centering
    \includegraphics[scale=0.3]{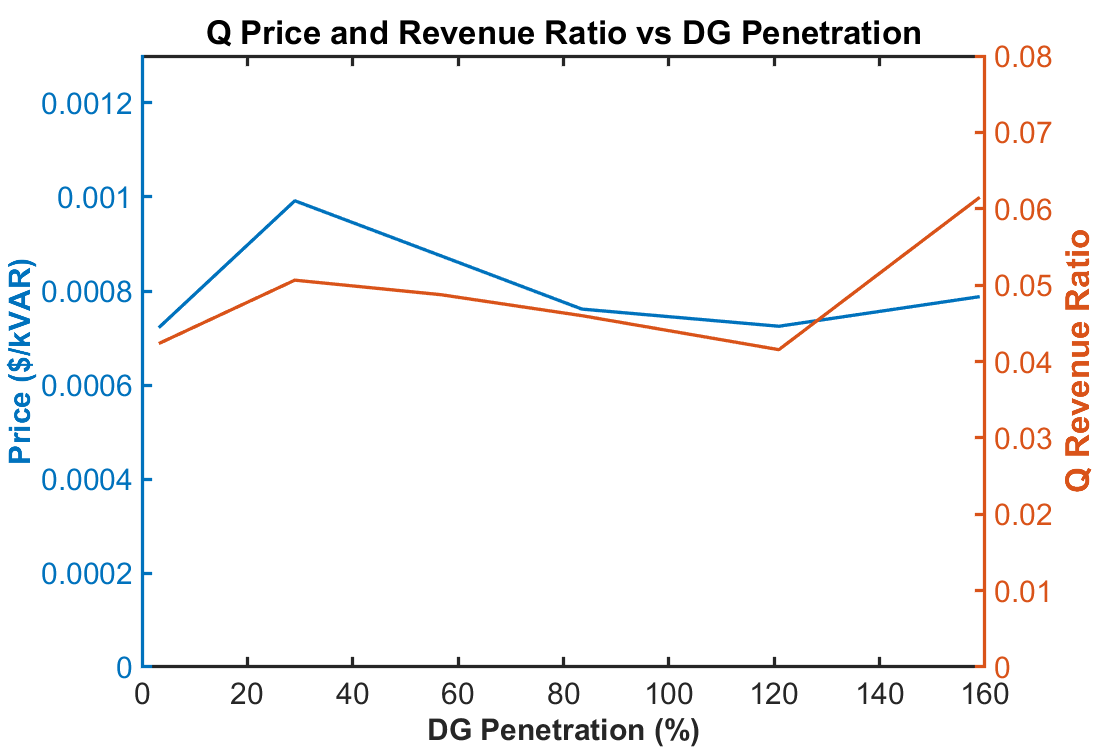}
    \caption{As the penetration of DGs increases, both the average price and Q-revenue ratio remain relatively stable. The price of reactive power is not strongly dependent on the PV penetration, making it more reliable for the transitioning grid. The revenue ratio shows that even as DGs are added, reactive power remains a reliable source of revenue. }
    \label{fig:penetration_graph}
\end{figure}

\begin{figure}[h!]
    \centering
    \includegraphics[scale=0.3]{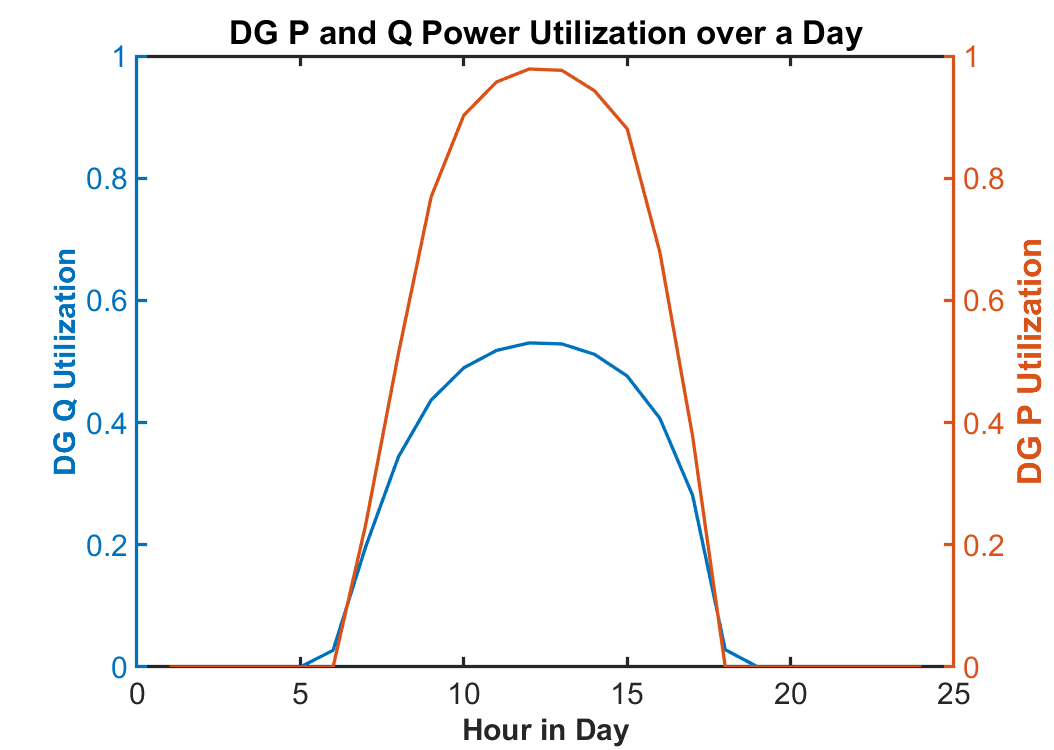}
    \caption{With a high DG penetration of 160\%, DGs nearly saturate the grid in the middle of the day. At the same time, DGs cover more than half of the reactive power load; much of the remaining load is covered by capacitor banks.}
    \label{fig:PQutil_day}
\end{figure}

Figure~\ref{fig:penetration_graph} shows the variation in $\bar{\mu}_j^Q$ and Q revenue, as a function of the DG penetration. It is clear that $\bar{\mu}_j^Q$ remains relatively constant even as the DG penetration increases from 5\% to 160\%. We also observe that the Q-revenue ratio remains relatively steady around 5\% meaning the importance of the reactive power market to DGs' revenue does not degrade. Figure~\ref{fig:PQutil_day} shows the DG-Q and DG-P utilization over the course of a day, for a DG penetration of 160\% (equivalent to a DG energy penetration of 35\%). In this scenario, when solar generation peaks around mid-day, almost 100\% of real power load can be met with local generation. This saturation in real power means that any additional DG capacity added to the grid will be used for meeting reactive power load. This shift, wherein DGs move from supplying primarily real power to reactive power, results in an increase in Q-revenue ratio. This is seen in Fig.~\ref{fig:penetration_graph}, by the Q-revenue ratio inflection point at high DG penetration. Notably, even at this high penetration, the reactive power price $\bar{\mu}_j^Q$ remains steady. Together, these results show that as more DGs are added, they can support the grid's real and reactive power needs without hurting the revenue potential of other DGs. This result may be due to power factor flexibility creating a dynamic supply of reactive power. As added DGs increase the overall supply of reactive power, prices can only drop so low before DGs choose to increase their power factors to sell in the real power market and vice versa. As the real and reactive power markets are cleared together, the two are interdependent such that a stable real power market promotes the stability of the reactive power market. In this way, our market maintains consistent payouts for reactive power, and can support the grid as it transitions to high DER penetration.

\begin{figure}[h!]
    \centering
    \includegraphics[scale=0.3]{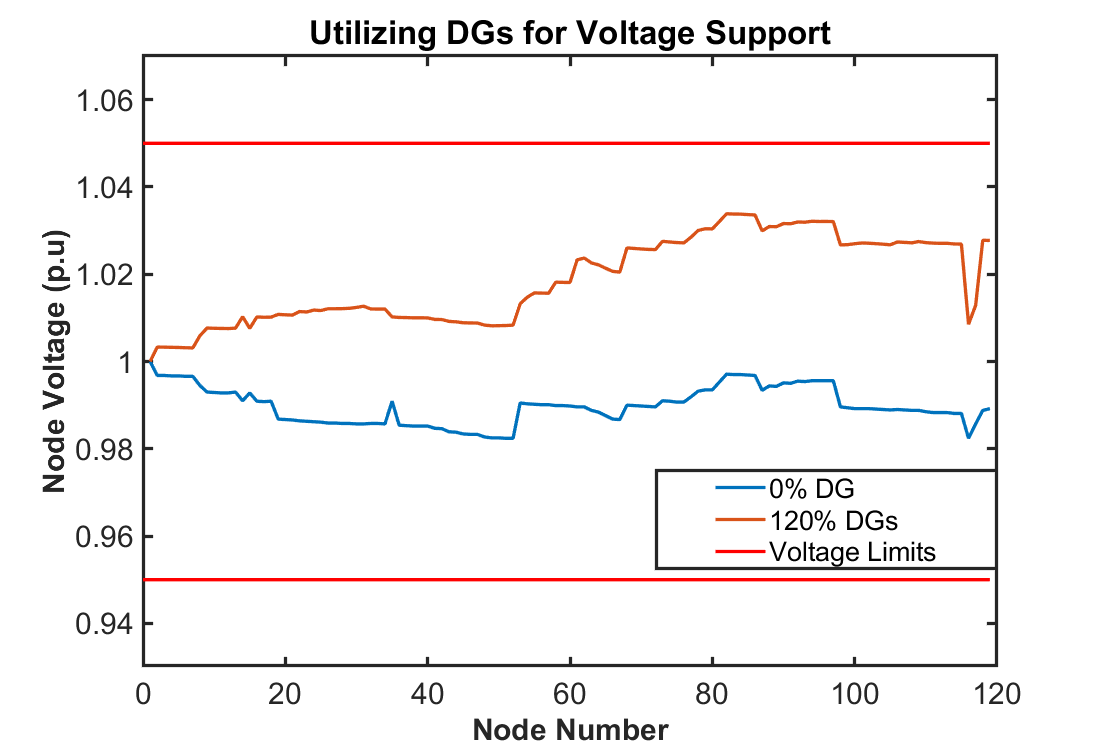}
    \caption{By adding DGs that provide Q at a PF of 0.9, the average node voltage of the grid increases by 2.2\% while staying within voltage limits. The DGs are able to provide necessary grid services by supporting voltages during periods of high demand.}
    \label{fig:voltage_support}
\end{figure}

Through these new revenue stream for DER owners, the reactive power market incentivizes DG resources to provide necessary services to the grid, including improved voltage characteristics.
Figure~\ref{fig:voltage_support} plots the average voltage profile over the network, when utilizing DGs within the PF of 0.9 to 1, for 1pm on June 27, 2021. We note the voltages across the feeder (at every phase and every node) remain within the specified voltage limits ($\pm5\%$) at all times, enforced by the voltage constraints in the CI-OPF formulation in \eqref{eq:CIOPF_Vlim1}-\eqref{eq:CIOPF_Vlim2}. However, the mean voltage for the 120\% penetration case is 2.2\% higher, attributed to the higher DG penetration and coordinated loss reduction enabled by the market. The voltage limits in red are given as a reference. Even with the limited PF range of just 0.9 to 1, the distributed coordination of DGs enabled by our PAC-based CI-OPF model leverages the flexibility of the DGs to greatly improve the average voltage profile, without any overvoltage events which can happen in networks with high DG penetration. In comparison, the 0\% DG Penetration case only receives voltage support from capacitor banks, which cannot respond quickly to varying grid conditions.

\begin{table}[]
\begin{tabular}{c|cccccc}
DG Penetration (\%)                                                                                           & 3.2                 & 29                     & 56.7                    & 83.5                      & 121                       & 159                       \\ \hline
\begin{tabular}[c]{@{}c@{}}Economic loss\\ \textit{(Normalized)}\end{tabular}                                        & 1.00                    & 0.949                     & 0.853                      & 0.776                      & 0.643                      & 0.564                      \\
\begin{tabular}[c]{@{}c@{}}Efficiency loss\\ \textit{(Normalized)}\end{tabular}                                      & 1.00                    & 0.966                     & 0.897                      & 0.851                      & 0.783                      & 0.783                      \\
\begin{tabular}[c]{@{}c@{}}Total loss\\ \textit{(Normalized)}\end{tabular}                                           & 1.00                    & 0.949                     & 0.854                      & 0.776                      & 0.644                      & 0.565                      \\ \hline
\multicolumn{1}{l|}{\begin{tabular}[c]{@{}l@{}}Loss Reduction compared \\ to 3.2\% DG Penetration\end{tabular}} & \multicolumn{1}{l}{ --} & \multicolumn{1}{l}{5.1\%} & \multicolumn{1}{l}{14.6\%} & \multicolumn{1}{l}{22.4\%} & \multicolumn{1}{l}{35.6\%} & \multicolumn{1}{l}{43.5\%}
\end{tabular}
\caption{ The objective function described in Equation \ref{eq:CIOPF_obj_econ} is split into an economic loss describing the cost to the grid and an efficiency loss for line losses (lower is better for both). Higher DER penetrations show significant reductions in objective function value indicating lower costs and more efficient grid operation. Values are normalized against the 3.2\% DG penetration scenario for each category.}
\label{tab:loss_table}
\end{table}

The objective function described in Eq.~\eqref{eq:CIOPF_obj_econ} is decomposed into an economic term reflecting the overall cost of supplying power on the grid and an efficiency term representing line losses. The optimal objective value for varying DG penetrations is shown in Table~\ref{tab:loss_table}. Each loss term and the total loss is normalized against the highest cost point (i.e. at lowest DG penetration of 3.2\%) for each category, and the percentage of loss reduction as compared to the lowest DG penetration is presented. Most notably, the proposed market structure is able to take advantage of the low cost DGs present in the network to reduce both economic costs and network losses. The economic costs drop by 43.6\% at the highest DG penetration while DG payout remain stable (as shown in Fig.~\ref{fig:penetration_graph}). As real and reactive power generation is moved closer to loads, line losses also drop considerably. We note that the efficiency term represents a small fraction of the total loss, reflecting traditional market structure where optimization is driven primarily by economic costs. The 43.5\% improvement in the social welfare (last row of the table) is strong evidence that the proposed market is able to use the DGs to reduce costs to customers and grid operators without degrading DG revenue, and while minimizing line losses.

\subsection{Market Utilization of Variable Power Factors}
In this section we investigate the market behaviour under varying PF settings of the DERs. We vary the PF operating range of all inverters, varying minimum PF from 0.6 to 0.95 by changing the minimum PF in the inverter P-Q model \eqref{eq:CIOPF_pf}. This is a realizable scenario as DGs like PV are being increasingly equipped with smart inverters with advanced power electronics making them capable of operating in such a range. Unlike existing methods, our market is able to utilize DGs at lower PFs to support the grid, while creating a strong economic incentive for their participation in the reactive power market.

\begin{figure}[h!]
    \centering
    \includegraphics[scale=0.3]{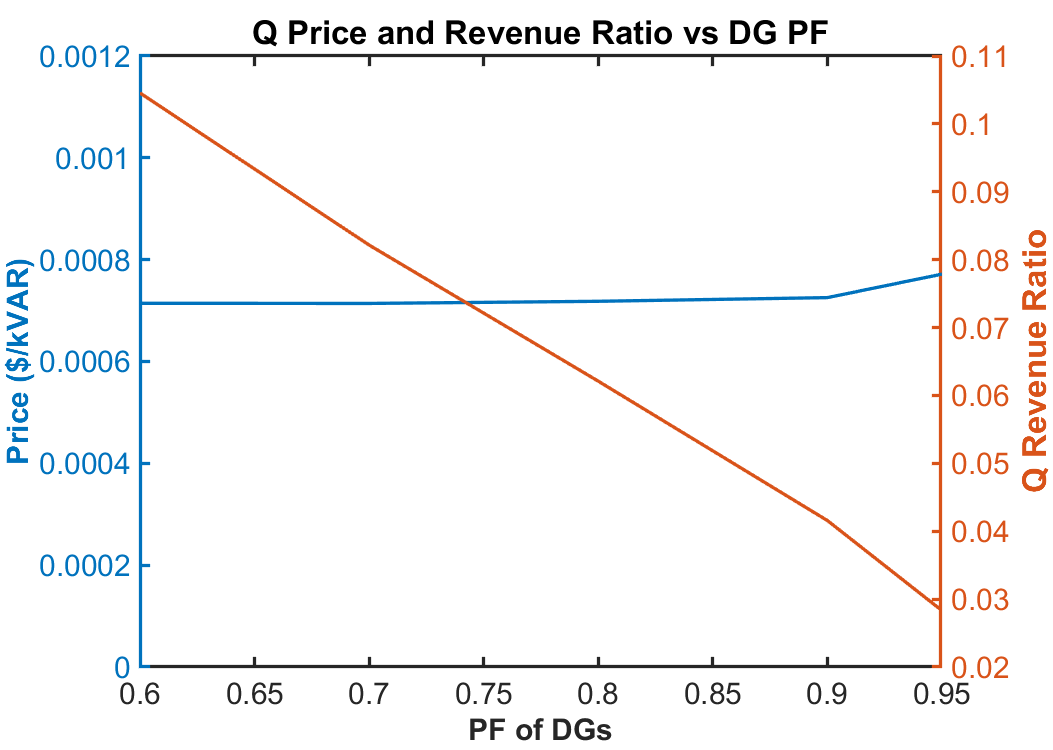}
    \caption{Enabling DGs to operate at low PFs increases the fraction of their revenue from reactive power. DGs that are permitted to have PFs as low as 0.6 sourced over 10\% of their revenue from our market's reactive power payments.}
    \label{fig:revenue_pf}
\end{figure}

Table~\ref{tab:DGtable} and Figure~\ref{fig:revenue_pf} show the corresponding results obtained for June 27, 2021, with a DG Penetration in \eqref{eq:penetration} of 120\%. The data in Table~\ref{tab:DGtable} illustrates how the DG-Q utilization varies with PF, and makes it apparent that by leveraging capabilities of inverters, upwards of 40\% of total reactive power loads can be served by the local DGs over a 24-hour period. As the inverter's operating range is increased, which corresponds to a change in minimum PF from 0.95 to 0.6, the DGs are capable of supplying more of the network's reactive power load, from 15\% to 44\%. Note that when operating at a PF of 0.6, the DG is supplying approximately 30\% more reactive power than real power, at no additional hardware costs to the generator. This is in contrast to synchronous generators which only operate at minimum PFs of 0.95, beyond which additional costs are incurred and efficiency degrades. This increased operating range at low costs means DGs can be an indispensable resource to the grid, when utilized optimally.

Figure~\ref{fig:revenue_pf} shows the variation in reactive power prices and Q-revenue ratio as a function of PF limit. We observed that $\bar{\mu}_j^Q$ (left axis of Figure \ref{fig:revenue_pf}) remained relatively stable, with only an 8\% decrease as minimum PF is reduced from 0.95 to 0.6. However, the Q-revenue ratio increases from 3\% to 10.5\%. This implies that DGs operating at lower PFs can anticipate up to 10\% of their total revenue to come from reactive power payments, which is a significant amount. In practice, DGs will chose a PF based on the relative demand of real and reactive power locally minimizing overall cost to the grid. Clearly, a retail market for reactive power, such as the one proposed in this paper, is essential to appropriately compensate DGs for their services. 

Under this scenario of variable PF over a wide range, the coordination of DGs by the market results in a notable improvement in the voltage profile. Table~\ref{tab:DGtable} outlines the mean voltage over the network and total network losses. As DG inverters are permitted to operate at lower PFs, the mean voltage across the network improves as the DG-Q utilization increases. Note that this comes at an expense of increased network losses: the objective function \eqref{eq:CIOPF_obj_econ} considers a trade-off between maximizing the social welfare and minimizing line losses, capture by parameter $\zeta$. At higher values of $\zeta$, the line losses are reduced more substantially.

Our results shows that our market will utilize the full flexibility in PF operating range offered by DGs through the variable payments, something not possible under current industry practice of fixed pricing.

\begin{table}[]
\begin{tabular}{c|c|c|c}
Power Factor & DG-Q Utilization & Mean Voltage (p.u.) & Network Losses (p.u.) \\ \hline
1.0          & 0.000            & 0.990               & 0.0344                \\
0.95         & 0.147            & 0.992               & 0.0353                \\
0.9          & 0.204            & 0.993               & 0.0361                \\
0.8          & 0.284            & 0.996               & 0.0383                \\
0.7          & 0.357            & 0.997               & 0.0399                \\
0.6          & 0.442            & 0.998               & 0.0406               
\end{tabular}
\caption{DGs are allowed to operate over a range of power factors. As the power factor deviates from 1.0 the Q utilization increases significantly, resulting in improved mean voltage with a modest increase in line losses.}
\label{tab:DGtable}
\end{table}

\subsection{Market Stability over a Week}
To better understand the market performance over a range of conditions, we simulate the market over the week of June 27, 2021 to July 3, 2021, using demand and LMP data from ISO-NE. This week is particularly interesting because it coincided with the end of a heat wave in New England and thus represents a larger than average variation in grid demand. For this simulation, a DG penetration of 120\% was used, corresponding to the addition of PV clusters at 20 nodes in the network, each with an inverter with a minimum PF of 0.9.

\begin{figure}[h!]
    \centering
    \includegraphics[scale=0.3]{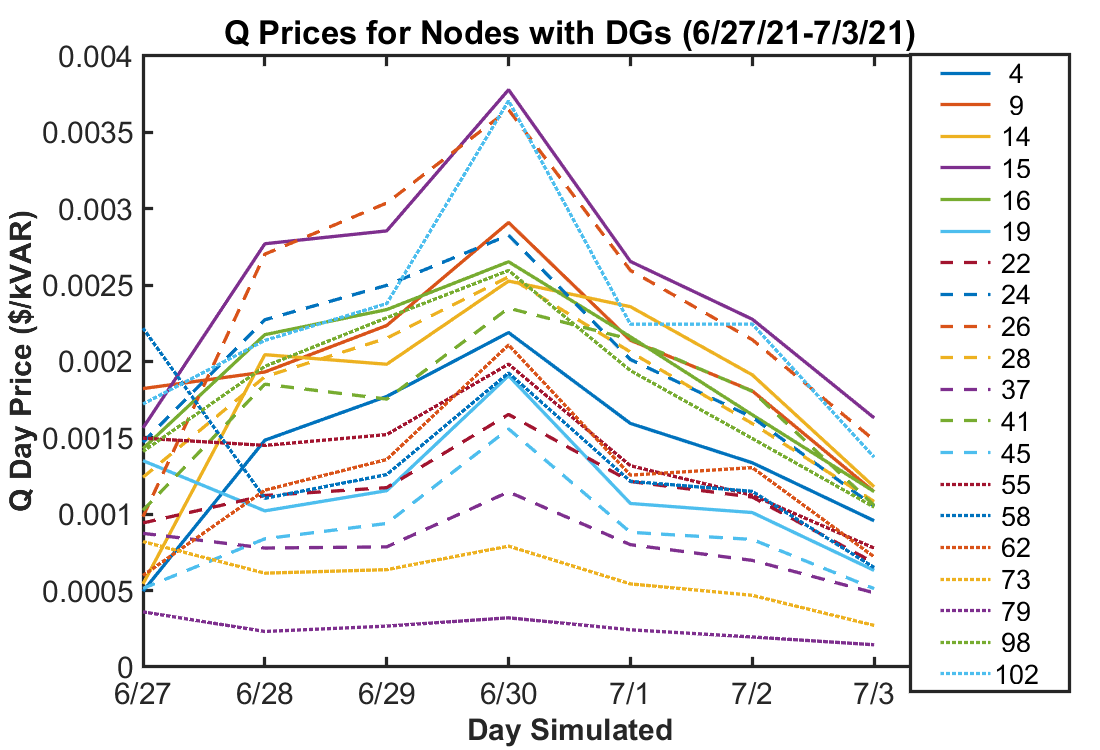}
    \caption{The daily prices of reactive power vary spatially across nodes and temporally over a week, but maintain the same order of magnitude. Prices remain relatively consistent, creating a stable and reliable revenue stream for DGs.}
    \label{fig:week_avg}
\end{figure}

Figure~\ref{fig:week_avg} shows the resulting $\bar{\mu}_j^Q$ values over the week, for all 20 DGs. The significant price differences between nodes speaks to the importance of d-LMPs due to their variability in both time and location. The temporal variation in price exhibited over the week further motivates the need for time-varying prices in order to reflect the varying grid conditions, especially with increased DER-penetration. These variations cause each node to fluctuate around an average price creating a signal for DG investments. For example, a DG at nodes 15 and 102 are likely to receive more revenue from reactive power than a DG at nodes 73 or 79.

\begin{figure}[h!]
    \centering
    \includegraphics[scale=0.32]{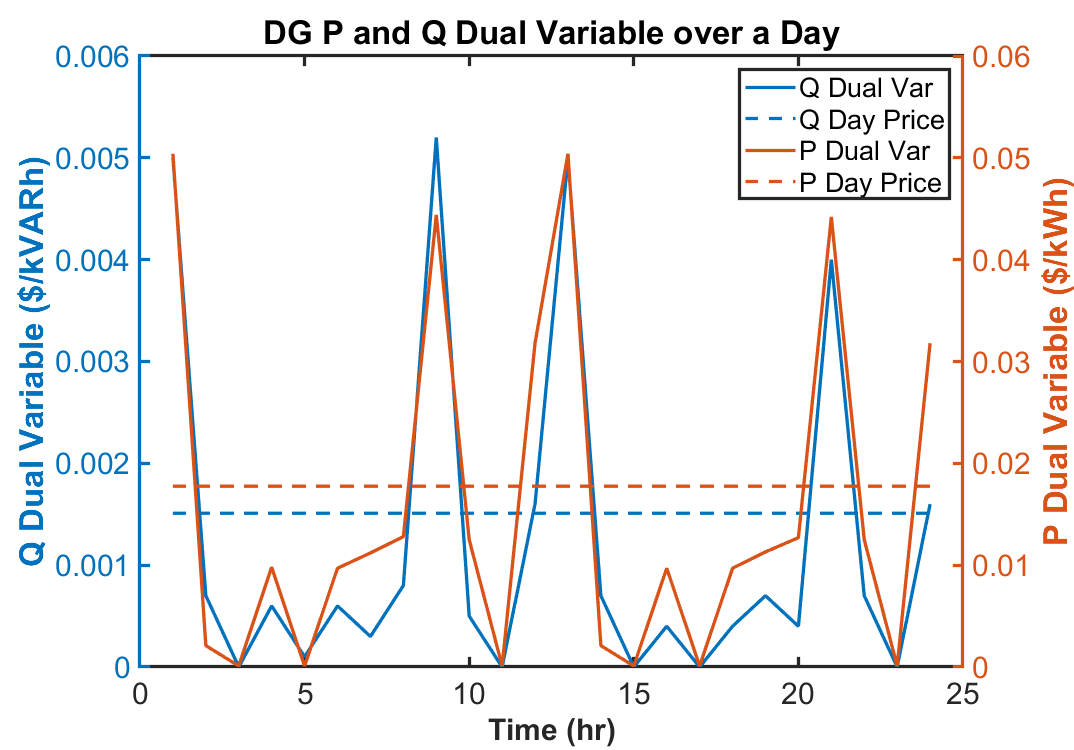}
    \caption{The hourly dual variables $\mu_{j,t}^Q$ and $\mu_{j,t}^Q$ (dual variables) shows significant volatility for the DG located at node 16 over the course of a single day (6/29/2021). In comparison, the proposed daily d-LMP, shown in dashed lines, are approximately 3.3 less volatile than the hourly price, while both result in the same payout (by construction). Note that the reactive power d-LMP is roughly 10\% of the real power d-LMP.}
    \label{fig:q_price}
\end{figure} 

Price volatility is a key concern for reactive power markets \cite{weber_simulation_1998,barquin_gil_reactive_2000,roozbehani_dynamic_2010}, and motivates the weighted average pricing model introduced in \eqref{eq:price_def}. The real and reactive power prices in Fig.~\ref{fig:week_avg} are calculated using the dual variables $\mu_{j,t}^P$ and $\mu_{j,t}^Q$ respectively, which are shown in Fig.~\ref{fig:q_price} for a single DG (node 16) over a 24-hr simulation. Compared to the stable daily prices, the underlying hour-to-hour price exhibits large price variations, something that is undesirable for a market. The average $c_v$ (coefficient of variation, a measure of price volatility) for $\bar{\mu}_j^Q$ is 0.42, approximately 3.3 times less volatile than the hourly $c_v$ of 1.39. The $\bar{\mu}_j^Q$ and $\bar{\mu}_j^P$ values are shown with dashed lines in Fig.~\ref{fig:q_price} for a particular DG and day. While the total payout to DERs is the same for the hourly duals $\mu_{j,t}^Q$ or using the daily d-LMPs $\bar{\mu}_j^Q$, the daily d-LMPs are significantly less volatile (similarly for real power). Averaging over each day, therefore, is chosen for our real and reactive power d-LMP, $\bar{\mu}_j^P$ and $\bar{\mu}_j^Q$, in order to contain the volatility and provide a reliable revenue stream for DERs. A choice of a larger period such a month or a year would suppress useful granularity in prices, and render an inefficient pricing mechanism. As shown in Fig.~\ref{fig:week_avg}, the daily reactive power prices as observed by DGs are not nearly at volatile as the corresponding hourly dual variable in Fig.~\ref{fig:q_price}. By managing the level of volatility the DER owners are exposed to, we instill confidence in the market while providing reliable revenue streams for DERs, thereby supporting the long-term development of DGs.

For a given DG, the relationship between real and reactive power hourly prices can be also observed in Fig.~\ref{fig:q_price}, shown here for node 16. Because the two are cleared simultaneously, real and reactive power hourly prices generally fluctuate together with reactive power at roughly 10\% of the real power price. This ratio of real to reactive power d-LMP at the distribution level mirrors the ratio of real to reactive power LMP at the wholesale level, where reactive power price is chosen to be one-tenth the real power LMP \cite{federal_energy_regulatory_commission_payment_2014}. 

\subsection{ Implementation Challenges}
The results of the previous section show how the growing penetration of DERs are capable of providing significant grid support through their adjustable PFs. Current market practices are inadequate to integrate and compensate small grid-edge devices, particularly for reactive power. Our proposed distribution level retail market achieves this by enabling DERs to provide grid services and be remunerated at a fair retail rate. To realize such a retail market structure requires significant effort in regulatory, policy, and technology investments and changes.

Our retail market is proposed to be overseen by a DSO, an entity which is not revenue driven and has regulatory and market oversight of DER participation, similar to the ISO at the transmission level. Rather than an aggregator, the DSO can be viewed as a proactive utility which optimally makes use of the DERs within the distribution network to maximize economic efficiency and minimize power losses. To implement such a DSO, the utility would require restructuring to focus on the integration, support, and utilization of distribution level assets. While there may be initial resistance to such a change, our simulation results show that utilities stand to gain from such a market, by better resource management including DERs. In one simulation scenario, tapping into the reactive power capabilities of DERs was able to aid utilities by provide over 45\% of reactive power needs. Our market is beneficial to both utilities and resource owners, with our results showing DER owners stand to receive over 10\% of their income from reactive power. In addition to the cooperation of existing utilities and load serving entities, the market incentives and process must encourage high rates of market enrollment by DER owners. The timeline and costs must be non-prohibitive. The question of who bears the cost of system upgrades for interconnection, additional metering (such as bidirectional meters), or upgrading inverters which follow updated operation standards (such as IEEE 1547) remains an open policy question for DER interconnection and wholesale market participation through aggregators (as required by FERC 2222). Similar questions will arise for the distribution-level retail market proposed herein.

Our market proposal introduces a price for reactive power, at which DG owners are compensated and loads are charged for their usage. As such, consumers must also become accustomed to this new pricing scheme for power. Currently, reactive power costs are represented on bills as delivery charges and generation service charges. In the future, reactive power rates could be stated separately with daily or weekly aggregate prices. This serves to inform customers of their usage patterns, and DER owners of how their devices are being used to support the grid. As with any new market derivative and changes to customer billing, the stability and fairness of pricing mechanisms must be thoroughly investigated. To this end, market stability has been explicitly examined in this paper, showing stable reactive power prices over a week of operation, and stable prices across DG penetration and minimum PFs. Further investigation of the proposed retail market is needed to ensure rate equity across different socioeconomic classes and fair access of electricity; we note that this is not unique to our reactive power market, but rather is a necessary feature of any electricity market design.

In order to implement the PAC algorithm to coordinate price negotiations between agents, the grid will require significant communication infrastructure. Fortunately, the shift towards smart meters, which report consumption metrics in real time across a bidirectional communication pathway, is well underway. A 2021 FERC report found about 60\% of meters in the US were smart meters \cite{ferc_smart_metering}. Implementation of the PAC algorithm and the subsequent reactive power market will likely require some additional hardware, software, and cybersecurity measures, but progress towards advanced metering infrastructure suggests that these additional implementations are becoming more and more practicable.

Finally, an important note on ancillary market design must be made. Reactive power control is the key tool used in voltage control and ensuring voltage stability of the grid. As previously mentioned, this paper deals with a reactive power energy market, and does not consider voltage control services. We leave to future work the consideration of how the proposed energy market can interface with ancillary markets at the distribution level. Our market takes the first step in enabling resource coordination and results in improved voltage profiles over the feeder by minimizing line losses. However, as DER penetration increases, resulting in more variability in load and generation, there may be a need for additional grid services like voltage control to more closely manage voltages. Corresponding market derivatives for reactive power may then be needed in an ancillary market framework. Further, coordination between distribution-level and wholesale-level ancillary markets must be made, and reactive power provisions ensured to maintain the stability of the interconnected grid.

\subsection{Key Observations} \label{sec:key_observ}

Through these simulations and analyses, the following key observations were made:

\begin{itemize}
     \item Our market allows DGs to contribute a significant amount of grid support by leveraging a wide range of PFs. In simulation, DGs accounted for 42\% of reactive power when allowed to operate at PFs as low as 0.6 (Table~\ref{tab:DGtable}).
    \item The reactive power market introduced and evaluated in this paper can incentivize DERs to participate in reactive power support with Q-revenue streams of up to 10\% of total revenue (Fig.~\ref{fig:revenue_pf}).
    \item The resource coordination enabled by our market structure uses DGs towards supporting network voltages, with average voltage boosted by up to 2.2\% as compared to zero DG penetration (Fig.~\ref{fig:voltage_support}, Table~\ref{tab:DGtable}).
    \item The limited price volatility of the daily Q d-LMPs (Fig.~\ref{fig:week_avg}) and the stable revenue streams under varying DG penetration (Fig.~\ref{fig:penetration_graph}) can instill confidence in the market, through the energy transition. This could encourage additional DG enrollment into the retail market, and also help drive investment decisions for additional DER adoption.
\end{itemize}

\section{Conclusion} \label{sec:conclusion}

In this paper, we have proposed a reactive power market to enable DG participation in the distribution grid. The market scheduling and clearing is based on PAC, a distributed optimization algorithm, which has been shown to converge in other works. Through a numerical case study of a modified IEEE-123 bus, we demonstrated that our market accommodates a large amount of DGs with a DG penetration anywhere between 5\% and 160\%, with the revenue stream remaining steady. We also showed that this market utilizes the full flexibility of DGs with the ability to operate over a range of PF, from 0.6 to 0.95 thereby meeting over 45\% of reactive power load, with the reactive power price per unit kVAR remaining the same. We also showed through a detailed assessment of the price variations across all DGs in the IEEE 123-bus over a week shows that the price variations are non-volatile; they primarily follow load fluctuations and otherwise exhibit little fluctuations.

The proposed market is structured so as to simultaneously clear real and reactive power with the resulting solution serving as the DER dispatch, and the duals providing a basis for real and reactive power pricing. It uses variable payments and spatially varying prices to better capture the value of DERs and more effectively control DERs through economic signaling. The reactive power market incentivizes DERs to enroll and participate in reactive power support by offering a reliable additional source of income through daily reactive power distributed locational marginal prices (Q d-LMPs). The results from our numerical study show that this revenue stream remains dependable even as DG penetration rapidly grows. By building confidence with stable prices, as shown in Section~\ref{sec:results}, additional DER enrollment may perhaps by encouraged into the retail market and help drive investment decisions for additional DER adoption. Our results appear promising in this regard and permit better DER integration and utilization to meet distribution grid objectives, and can enable deeper DER penetration.

\normalem

\bibliography{references}{}


\end{document}